\documentclass{amsart}
\usepackage[T1]{fontenc}
\usepackage{lmodern}
\usepackage[a4paper,margin=29mm]{geometry}
\usepackage{amsmath,amssymb,amsthm,mathtools}
\usepackage{booktabs,array}
\usepackage{capt-of}
\usepackage{enumitem}
\usepackage{microtype}
\usepackage{xurl}
\usepackage[hidelinks]{hyperref}
\hypersetup{pdftitle={The optimal exponent in the Babai--Goodman--Pyber theorem},
  pdfauthor={Shaowu Zheng},
  pdfkeywords={Finite groups, soluble subgroups, soluble radical, finite simple groups, Babai--Goodman--Pyber theorem}}
\newtheorem{theorem}{Theorem}[section]
\newtheorem{lemma}[theorem]{Lemma}
\newtheorem{proposition}[theorem]{Proposition}

\newtheorem*{theoremA}{Theorem A}
\newtheorem*{theoremB}{Theorem B}
\theoremstyle{definition}

\theoremstyle{remark}
\newtheorem{remark}[theorem]{Remark}
\DeclareMathOperator{\Aut}{Aut}
\DeclareMathOperator{\Out}{Out}
\DeclareMathOperator{\Rad}{Rad}
\DeclareMathOperator{\Soc}{Soc}
\DeclareMathOperator{\GL}{GL}
\DeclareMathOperator{\SL}{SL}
\DeclareMathOperator{\PSL}{PSL}
\DeclareMathOperator{\PGL}{PGL}

\DeclareMathOperator{\SU}{SU}
\DeclareMathOperator{\PSU}{PSU}
\DeclareMathOperator{\Sp}{Sp}
\DeclareMathOperator{\PSp}{PSp}
\DeclareMathOperator{\ord}{ord}

\newcommand{\PO}{\mathrm{P}\Omega}
\newcommand{\F}{\mathbb F}
\allowdisplaybreaks[2]
\numberwithin{equation}{section}
\title{The optimal exponent in the Babai--Goodman--Pyber theorem}
\author{Shaowu Zheng}
\address{Department of Mathematics, University of Auckland, Private Bag 92019, Auckland 1142, New Zealand}
\email{shaowu.zheng@auckland.ac.nz}
\date{November 2026}

\begin{document}

\begin{abstract}
Babai, Goodman and Pyber proved that there exists an absolute
constant $c$ such that
$[X:\operatorname{Rad}(X)]\le [X:H]^c$
for every finite group $X$ and every soluble subgroup $H\le X$,
where $\operatorname{Rad}(X)$ denotes the largest soluble normal
subgroup of $X$.
We prove that the least possible value of $c$ is
\[
b_0=1+\log_9(168\cdot 24^{1/3})\approx3.8141.
\]
For nonsoluble $X$, we obtain the stronger inequality
$[X:\operatorname{Rad}(X)]\le 24^{-1/3}[X:H]^{b_0}$,
and the coefficient $24^{-1/3}$ is optimal for this exponent.
We reduce the problem to bounds on soluble subgroups of
automorphism groups of finite simple groups.
Examples based on $\operatorname{PSL}_2(8)$
establish sharpness.
\end{abstract}
\maketitle

\section{Introduction and main results}\label{sec:intro}
Let $X$ be a finite group and let $H\le X$ be soluble. Write
$\Rad(X)$ for the largest soluble normal subgroup of $X$.
Babai, Goodman and Pyber \cite[Theorem~2.13]{BGP} proved that there exists an absolute
constant $c$ such that
\[
 [X:\Rad(X)]\le [X:H]^c.
\]
They conjectured that one may take $c=7$
\cite[Conjecture~6.6]{BGP}.
A stronger conjecture, attributed to Babai, Goodman and Pyber
\cite[Problem~17.41(a)]{Kourovka}, states that, when $\Rad(X)=1$,
seven conjugates of $H$ have trivial intersection.
This would imply that $c=7$ suffices.
Vdovin later conjectured that five conjugates suffice
under the same hypothesis
\cite[Problem~17.41(b)]{Kourovka}.
For finite primitive permutation groups with soluble point
stabilisers, Burness \cite[Theorem~1]{BurnessPrimitive}
proved that there exist five conjugates of a point
stabiliser whose intersection is trivial.
For prime-power index $m=[X:H]=p^a$, Bastos and Schneider
\cite[Theorem~4]{BastosSchneider} proved that $[X:\Rad(X)]\le m^5$, and
improved this to $[X:\Rad(X)]\le m^4$ when $p\ne13$.

We determine the least possible value of $c$ for all finite
groups $X$ and soluble subgroups $H\le X$.
Let
\[
\delta=24^{1/3},\qquad
r_0=\log_9(168\delta),\qquad
b_0=1+r_0\approx3.8141.
\]

\begin{theoremA}
For every finite group $X$ and every soluble subgroup $H\le X$,
\[
 [X:\Rad(X)]\le [X:H]^{b_0}.
\]
If $X$ is nonsoluble, then the stronger inequality
$
 [X:\Rad(X)]\le \delta^{-1}[X:H]^{b_0}
$
holds. The exponent $b_0$ is the least possible in the first inequality.
For this exponent, $\delta^{-1}$ is the least possible coefficient
in the second inequality.
\end{theoremA}
\begin{theoremB}
Let $r\ge0$. The following are equivalent.
\begin{enumerate}[label=\textup{(\roman*)}]
\item For every nonabelian finite simple group $T$ and every soluble
$B\le\Aut(T)$, $\delta|B|\le[T:B\cap T]^r$.
\item For every finite group $X$ and every soluble $H\le X$,
$[X:\Rad(X)]\le[X:H]^{r+1}$.
\end{enumerate}
Moreover, if \textup{(i)} holds, then
$[X:\Rad(X)]\le\delta^{-1}[X:H]^{r+1}$
for every nonsoluble finite group $X$
and every soluble subgroup $H\le X$.
\end{theoremB}

By Theorem~B with $r=r_0$, the inequalities in Theorem~A follow
once we show that
$\delta|B|\le[T:B\cap T]^{r_0}$ for every nonabelian finite simple
group $T$ and every soluble subgroup $B\le\Aut(T)$.
The proof of the implication (i)$\Rightarrow$(ii) in Theorem~B uses
the nonabelian socle of $X/\Rad(X)$ and Dixon's bound for soluble
permutation groups \cite[Theorem~3]{Dixon}. The proof of the converse
uses iterated imprimitive wreath products with $\operatorname{Sym}(4)$.
The sharpness of the exponent and the coefficient
in Theorem~A follows from Theorem~B and Proposition~\ref{prop:lower}.

To establish this inequality,
we prove a uniform estimate for nonabelian
finite simple groups (Theorem~\ref{thm:uniform}) by
induction on the order of the simple group. We use the classification of finite simple
groups to treat the alternating, sporadic and Lie-type families. For
the sporadic groups and the Tits group we use Burness's result
\cite{Burness}. 
We use Mann's results~\cite{Mann} to determine
the maximum orders of soluble subgroups of
$\PSL_n(q)$ for $n\ge3$.
Minimum-index bounds~\cite{GMPS} treat $\PSU_3(q)$, $\PSU_5(q)$
and $\PSp_4(q)$, including $\PSU_4(2)\cong\PSp_4(3)$.
For the remaining classical groups and the split exceptional groups
$E_6(q)$, $E_7(q)$ and $E_8(q)$, in characteristic $p$, we
consider separately soluble subgroups $U$ with $O_p(U)\ne1$ and those
with $O_p(U)=1$. In the first case, the Borel--Tits theorem
\cite{GLS3,BW} places $U$ in a parabolic subgroup, which allows us to
apply induction to the nonabelian composition factors of its Levi
quotient. In the second case, we use bounds for
soluble linear groups and torus normalisers
\cite{Palfy,Wolf,Vdovin}.
The remaining seven exceptional families are treated using lower
bounds for the index of a proper subgroup.

Appendix~\ref{app:verification} gives the remaining
arithmetic checks for $\PSL_2(q)$ and the enumeration
of the remaining classical parabolic cases.

The paper is organised as follows.
Section~\ref{sec:ltg} proves the local-to-global equivalence and
sharpness. Section~\ref{sec:principles} gives the uniform reduction.
Sections~\ref{sec:small}--\ref{sec:exceptional} establish the estimates
for the simple-group families. Finally, Section~\ref{sec:proofs}
completes the proof of Theorem~A.

\section{The local-to-global equivalence and sharpness}\label{sec:ltg}

\begin{proof}[Proof of Theorem~B]

Assume (i). Let $X$ be a finite group and let $H\le X$
be soluble. Write $R=\Rad(X)$, $G=X/R$ and $J=HR/R$.
Then $J$ is soluble and $\Rad(G)=1$.
If $G=1$, then $[X:R]=1$, so (ii) holds.
Suppose that $G\ne1$.
Since $H\le HR$, we have $[G:J]=[X:HR]\mid[X:H]$.
Since $[X:R]=|G|=[G:J]|J|$, it suffices to prove
$|J|\le\delta^{-1}[G:J]^r$.

To apply (i), we consider the action of $J$ on the
simple factors of the socle of $G$.

\medskip
\noindent\emph{Step 1. The action on the simple factors.}

Since $\Rad(G)=1$, the socle of $G$ has the form
$S=\Soc(G)=T_1\times\cdots\times T_k$, where $k\ge1$ and
all $T_i$ are nonabelian simple groups.
Since $S\trianglelefteq G$, conjugation by $J$ induces
automorphisms of $S$. As $T_1,\ldots,T_k$ are the
minimal normal subgroups of $S$, these automorphisms permute
the factors.
Let $K$ be the kernel of the action of $J$ on
$\{T_1,\ldots,T_k\}$, and let $P\le\operatorname{Sym}(k)$
be the permutation group induced by $J$.
Then $P\cong J/K$ is soluble and $|J|=|K||P|$.
We bound $|K|$ using (i), and $|P|$ using Dixon's bound.

\medskip
\noindent\emph{Step 2. Bounding $|K|$.}

Since $K$ normalises each $T_i$, conjugation gives a
homomorphism $K\to\Aut(T_i)$.
Let $B_i$ be its image, and identify $T_i$ with its inner
automorphism group in $\Aut(T_i)$.
Write $W_i=B_i\cap T_i$. Since $B_i$ is a homomorphic image
of the soluble group $K$, it is soluble.

We first show that $K$ acts faithfully on $S$.
Since $S\trianglelefteq G$, we have $C_G(S)\trianglelefteq G$.
If $C_G(S)\ne1$, it contains a minimal normal subgroup $N$ of $G$.
By the definition of the socle, $N\le S$, and hence
$N\le S\cap C_G(S)=Z(S)=1$, a contradiction.
Thus $C_G(S)=1$.
The homomorphisms $K\to B_i$ induce a homomorphism
$K\to\prod_{i=1}^k B_i$.
Since its kernel is $C_K(S)\le C_G(S)=1$,
this homomorphism is injective. Hence
$|K|\le\prod_{i=1}^k|B_i|$.
Applying (i) to each $B_i$, we obtain
\[
|K|\le\prod_{i=1}^k|B_i|
\le\delta^{-k}
\left(\prod_{i=1}^k[T_i:W_i]\right)^r.
\]

We now bound the product of indices in terms of $[G:J]$.
Since every $T_i$ is normal in $S$, every element of $S$
normalises each $T_i$. By the definition of $K$,
we have $J\cap S\le K$.
For $x=(x_1,\ldots,x_k)\in J\cap S$, conjugation by $x$
on $T_i$ is the inner automorphism induced by $x_i$.
Since $x\in K$, this automorphism belongs to $B_i$,
so $x_i\in W_i$. Hence $J\cap S\le\prod_{i=1}^k W_i$.
Consequently,
\[
\prod_{i=1}^k[T_i:W_i]
=[S:\prod_{i=1}^k W_i]
\mid[S:J\cap S]=[SJ:J]\mid[G:J].
\]
Since $r\ge0$, we conclude that
$|K|\le\delta^{-k}[G:J]^r$.

\medskip
\noindent\emph{Step 3. Completing the bound.}

Since $P$ is a soluble permutation group of degree $k$,
Dixon's bound \cite[Theorem~3]{Dixon} gives
$|P|\le24^{(k-1)/3}=\delta^{k-1}$.
Combining this with the bound for $K$, we obtain
$|J|=|K||P|\le\delta^{-1}[G:J]^r$.
Therefore
\[
[X:R]=|G|=[G:J]|J|
\le\delta^{-1}[G:J]^{r+1}
\le\delta^{-1}[X:H]^{r+1}.
\]
This proves the final assertion and, since $\delta>1$, also (ii).

\medskip

Conversely, assume (ii). Let $T$ be a nonabelian finite
simple group and let $B\le\Aut(T)$ be soluble.
Identifying $T$ with its inner automorphism group,
write $h=[T:B\cap T]$. We aim to prove
$\delta|B|\le h^r$.
Let $A=TB\le\Aut(T)$. Since $T\trianglelefteq\Aut(T)$,
we have $[A:B]=h$ and $|A|=h|B|$.
Thus it suffices to prove $\delta|A|\le h^{r+1}$.

To prove this, we construct finite groups $X$ with
$\Rad(X)=1$ and soluble subgroups $H\le X$,
and then apply (ii).
Let $n\ge1$ and let $Q\le\operatorname{Sym}(n)$ be soluble.
Let $X=A^n\rtimes Q$ and $H=B^n\rtimes Q$,
where $Q$ permutes the coordinates.
Then $H$ is soluble and $[X:H]=h^n$.

We first prove that $\Rad(X)=1$.
Let $M=T^n\trianglelefteq X$.
If $N\trianglelefteq X$ is soluble, then $N\cap M$
is a soluble normal subgroup of $M$.
Since $M$ is a direct product of nonabelian simple
groups, we have $N\cap M=1$.
Since both $N$ and $M$ are normal in $X$, we have
$[N,M]\le N\cap M=1$. Hence $N\le C_X(M)$.
Since $Z(T)=1$, we have
$C_{\Aut(T)}(\operatorname{Inn}(T))=1$,
and hence $C_A(T)=1$.
Every element of $C_X(M)$ normalises each simple factor
of $M$. Since $A^n$ normalises each of these factors and $Q$
permutes them faithfully, it follows that $C_X(M)\le A^n$.
Consequently,
$C_X(M)=C_{A^n}(M)=C_A(T)^n=1$.
Thus $N=1$, proving that $\Rad(X)=1$.
Applying (ii) to $X$ and $H$, we obtain
\[
|A|^n|Q|=|X|
\le[X:H]^{r+1}=h^{n(r+1)}.
\]
Taking $n$-th roots gives
$|A||Q|^{1/n}\le h^{r+1}$.
It remains to choose soluble permutation groups $Q$
of degree $n$ such that $|Q|^{1/n}\to\delta$.
Let $Q_0=1$, acting on one point, and define
$Q_{s+1}=Q_s\wr\operatorname{Sym}(4)$ in its
imprimitive action.
Each $Q_s$ is soluble and has degree $n_s=4^s$.
Since $|Q_{s+1}|=24|Q_s|^4$, induction gives
$|Q_s|=24^{(n_s-1)/3}=\delta^{n_s-1}$.
Applying the preceding inequality with $Q=Q_s$,
we obtain
$|A|\delta^{1-1/n_s}\le h^{r+1}$.
Since $n_s\to\infty$, it follows that
$\delta|A|\le h^{r+1}$.
As $|A|=h|B|$, we conclude that
$\delta|B|\le h^r$, proving (i).

\end{proof}

\begin{proposition}\label{prop:lower}
Let $T=\PSL_2(8)$ and let
$B=\mathrm{A\Gamma L}_1(8)\le\Aut(T)$
be a point stabiliser on $\mathbb P^1(\F_8)$.
For $X=TB$ and $H=B$, we have $\Rad(X)=1$ and
$
[X:\Rad(X)]=\delta^{-1}[X:H]^{b_0}.
$
\end{proposition}

\begin{proof}
We have $|B|=8\cdot7\cdot3=168$,
$|B\cap T|=56$ and $|T|=504$.
Hence $[X:H]=[T:B\cap T]=9$.
Since $X$ is almost simple, $\Rad(X)=1$.
Moreover, $|X|=9\cdot168$ and
$9^{b_0}=9\cdot168\delta$. Thus
\[
[X:\Rad(X)]=|X|
=\delta^{-1}9^{b_0}
=\delta^{-1}[X:H]^{b_0}.
\]
\end{proof}

\section{Uniform reduction for finite simple groups}
\label{sec:principles}

By Theorem~B, it suffices to prove
$\delta|B|\le[T:B\cap T]^{r_0}$ for every
nonabelian finite simple group $T$ and every
soluble subgroup $B\le\Aut(T)$.

\begin{lemma}\label{lem:crude}
Let $T$ be a nonabelian finite simple group,
and let $s>0$ and $k>1$. Suppose that
\begin{equation}\label{eq:crude}
\delta|\Out(T)|\,|T|<k^{1+s}.
\end{equation}
Then every soluble subgroup $B\le\Aut(T)$ with
$[T:B\cap T]\ge k$ satisfies
$
\delta|B|<[T:B\cap T]^s.
$
\end{lemma}

\begin{proof}
Write $h=[T:B\cap T]$.
Since $B/(B\cap T)\cong BT/T\le\Out(T)$, we have
\[
\delta|B|
\le\frac{\delta|\Out(T)|\,|T|}{h}
<\frac{k^{1+s}}{h}
\le h^s.
\]
\end{proof}

For a finite group $X$, define
$\sigma(X)=\max\{|K|:K\le X,\ K\text{ soluble}\}$.
For the rest of the paper, fix $s=39/14<r_0$.
For a nonabelian finite simple group $T$, consider the inequalities
\begin{align}
\sigma(T)^{1+s}&<|T|^s,
\label{eq:inner-bound}\\
\delta|\Out(T)|\sigma(T)^{1+s}&<|T|^s.
\label{eq:outer-bound}
\end{align}

Suppose that $T$ satisfies \eqref{eq:outer-bound}, and let
$B\le\Aut(T)$ be soluble.
Since $|B\cap T|\le\sigma(T)$, we have
$[T:B\cap T]\ge |T|/\sigma(T)>1$.
Taking $k=|T|/\sigma(T)$, condition
\eqref{eq:crude} is equivalent to
\eqref{eq:outer-bound}.
Hence Lemma~\ref{lem:crude} and $s<r_0$ give
$\delta|B|<[T:B\cap T]^s<[T:B\cap T]^{r_0}$.

\begin{theorem}\label{thm:uniform}
Every nonabelian finite simple group $T$ satisfies
\eqref{eq:inner-bound}. Moreover,
\eqref{eq:outer-bound} fails if and only if
\[
T\cong\PSL_2(8),\quad
\PSL_2(16),\quad\text{or}\quad\PSL_3(3).
\]
\end{theorem}

The three exceptional groups will be treated separately.
Some results below assume that
\eqref{eq:inner-bound} holds for all
nonabelian simple groups of smaller order.
This assumption will be justified by induction
on $|T|$ in Section~\ref{sec:proofs}.

For a finite group $X$, let $c(X)$ be the product of the orders
of its abelian composition factors, counted with multiplicity.
\begin{lemma}\label{lem:sigma}
Let $X$ be a finite group, and let $S_1,\ldots,S_m$
be its nonabelian composition factors, counted with
multiplicity. 
Then
$
\sigma(X)\le c(X)\prod_{i=1}^m\sigma(S_i).
$
\end{lemma}

\begin{proof}
If $N\unlhd X$, then for every soluble $K\le X$,
$
|K|=|K\cap N|\,|KN/N|
\le\sigma(N)\sigma(X/N).
$
Hence $\sigma(X)\le\sigma(N)\sigma(X/N)$.

Let $1=X_0\lhd X_1\lhd\cdots\lhd X_t=X$ be a composition series
of $X$. Applying the preceding inequality to $X_{i-1}\unlhd X_i$
for $i=1,\ldots,t$, we obtain
$\sigma(X)\le\prod_{i=1}^t\sigma(X_i/X_{i-1})$.
Here $\sigma(X_i/X_{i-1})=|X_i/X_{i-1}|$ when $X_i/X_{i-1}$ is
abelian. Since the product of the orders of the abelian factors
is $c(X)$, the result follows.
\end{proof}

\section{\texorpdfstring{$\PSL_2(q)$}{PSL(2,q)}, alternating and sporadic groups}\label{sec:small}
In this section, we establish \eqref{eq:outer-bound}
for the groups $\PSL_2(q)$ with $q\notin\{8,16\}$,
the alternating and sporadic simple groups, and
the Tits group. 
For $\PSL_2(8)$ and $\PSL_2(16)$,
we prove
$\delta|B|\le[T:B\cap T]^{r_0}$
for every soluble $B\le\Aut(T)$.

\begin{lemma}\label{lem:pgl2}
Let $1\ne H\le\PGL_2(q)$ be soluble, where $q=p^f$.
Then one of the following holds:
\begin{enumerate}[label=\textup{(\arabic*)}]
\item $H$ fixes a point of $\mathbb P^1(\F_q)$;
\item $H$ is contained in a torus normaliser of order
$2(q-1)$ or $2(q+1)$;
\item $p$ is odd and $|H|\mid24$.
\end{enumerate}
Moreover, every torus normaliser $N\le\PGL_2(q)$ of order
$2(q\pm1)$ satisfies $|N\cap\PSL_2(q)|=2(q\pm1)/\gcd(2,q-1)$.
\end{lemma}

\begin{proof}
Let $E$ be a minimal nontrivial normal subgroup of $H$.
Since $H$ is soluble, $E$ is elementary abelian
of order $r^a$ for some prime $r$.

\smallskip
\noindent
\textit{Case 1: $r=p$.}
Every nonidentity element of $E$ is unipotent and
has a unique fixed point in $\mathbb P^1(\F_q)$.
Since $E$ is abelian, all these elements have the
same fixed point. As $E\lhd H$, this point is
also fixed by $H$, giving \textup{(1)}.

\smallskip
\noindent
\textit{Case 2: $r\ne p$ is odd.}
Choose $1\ne x\in E$.
Then $x$ has two distinct eigenlines over
$\overline{\F}_q$, and its projective centraliser
is a torus.
Since $E$ is abelian, it lies in this centraliser,
so all nonidentity elements of $E$ have the
same pair of eigenlines.
As $E\lhd H$, conjugation by $H$ permutes
the nonidentity elements of $E$.
Hence $H$ preserves the pair of eigenlines of $x$.
This pair is Frobenius-stable since $x\in\PGL_2(q)$.
Its stabiliser in $\PGL_2(q)$ is a torus normaliser
of order $2(q-1)$ or $2(q+1)$.
Thus \textup{(2)} holds.

\smallskip
\noindent
\textit{Case 3: $r=2\ne p$.}
Choose an involution $x\in E$.
After conjugation over $\overline{\F}_q$,
we may assume that $x:z\mapsto-z$.
Its projective centraliser consists of the maps
$z\mapsto az$ and $z\mapsto a/z$.
Let $D$ be the torus consisting of the maps
$z\mapsto az$, which has index two in this
centraliser.
Since $D\cong\overline{\F}_q^\times$, every finite subgroup
of $D$ is cyclic. As $E$ is elementary abelian and centralises
$x$, it follows that $E\cap D=\langle x\rangle$.
Hence $|E|\le4$.

If $|E|=2$, then $x$ is central in $H$.
Thus $H$ preserves the eigenline pair of $x$,
giving \textup{(2)}.

If $|E|=4$, choose $y\in E\setminus D$.
Then $E=\langle x,y\rangle$, where
$y:z\mapsto b/z$ for some $b\ne0$.
A direct calculation shows that
$C_{\PGL_2(\overline{\F}_q)}(E)=E$.
Since $E\lhd H$, conjugation gives a homomorphism
from $H$ to $\Aut(E)$ with kernel $C_H(E)=E$.
Hence $H/E$ embeds into
$\Aut(E)\cong\operatorname{Sym}(3)$.
As $|E|=4$ and $|\Aut(E)|=6$,
we obtain $|H|\mid24$, giving \textup{(3)}.

It remains to prove the final assertion.
For even $q$, we have $\PSL_2(q)=\PGL_2(q)$, so there is
nothing to prove. For odd $q$, the group $\PSL_2(q)$ is the
kernel of the determinant square-class map
\[
\PGL_2(q)\longrightarrow
\F_q^\times/(\F_q^\times)^2.
\]
Both split and nonsplit tori contain elements
of nonsquare determinant. This is clear for
split tori, while for nonsplit tori it follows
from the surjectivity of the norm map
$\F_{q^2}^\times\to\F_q^\times$.
Consequently, each torus normaliser $N$ satisfies
$
[N:N\cap\PSL_2(q)]=2.
$
\end{proof}

\begin{proposition}\label{prop:rankone}
Let $T=\PSL_2(q)$, where $q=p^f\ge4$ and $q\ne5$.
For every soluble $B\le\Aut(T)$, write
$h=[T:B\cap T]$, and let
$r_{16}=\log_{17}(960\delta)<r_0$.
Then the following hold:
\begin{enumerate}[label=\textup{(\roman*)}]
\item If $q\notin\{8,16\}$, then
\eqref{eq:outer-bound} holds.
\item If $q=8$, then $\delta|B|\le h^{r_0}$.
\item If $q=16$, then $\delta|B|\le h^{r_{16}}$.
\end{enumerate}
In \textup{(ii)} and \textup{(iii)}, equality holds
if and only if
$B=\Aut(T)_\alpha$ for some
$\alpha\in\mathbb P^1(\F_q)$.
\end{proposition}

\begin{proof}
Write $d=\gcd(2,q-1)$. We have
$\Aut(T)=\mathrm{P\Gamma L}_2(q)$,
$|\Out(T)|=df$ and
$|\Out(T)||T|=fq(q^2-1)$.
For $q\notin\{8,16\}$, the bounds
\[
\sigma(T)\le\frac{|T|}{k}
\qquad\text{and}\qquad
\delta fq(q^2-1)<k^{1+s}
\]
together imply \eqref{eq:outer-bound}.
We take $k=q+1$ for $q\ne7$ and $k=7$
for $q=7$. We first establish the bound
$\sigma(T)\le |T|/k$ by showing that
$[T:U]\ge k$ for every soluble $U\le T$.

This holds for $U=1$, so suppose $U\ne1$
and apply Lemma~\ref{lem:pgl2} to $U\le\PGL_2(q)$.
In case~\textup{(1)}, $U$ lies in a point stabiliser
of $T$. Since such a stabiliser has order $q(q-1)/d$,
we have $[T:U]\ge q+1$.
In case~\textup{(2)}, $U\le T\cap N$ for some torus
normaliser $N\le\PGL_2(q)$ of order $2(q\pm1)$.
By the final assertion of
Lemma~\ref{lem:pgl2}, we have
$|T\cap N|=2(q\pm1)/d$.
Hence $[T:U]\ge q(q-1)/2\ge q+1$.
In case~\textup{(3)}, $q$ is odd and $|U|\le24$,
so $[T:U]\ge q(q^2-1)/48$.
The bound $q(q^2-1)/48$ is at least $q+1$
for $q\ge9$ and equals $7$ for $q=7$.

\smallskip
\noindent
\textit{(i) $q\notin\{8,16\}$.}
By the index bound just proved, it remains
to verify \eqref{eq:crude}, namely,
$\delta fq(q^2-1)<k^{1+s}$.
For $q\ne7$, this is equivalent to
$\delta fq(q-1)<(q+1)^s$.
For $q\ge64$, we have $f\le\sqrt q$.
Since $s-5/2=2/7$ and $3^7<64^2\le q^2$,
we also have $3<q^{s-5/2}$.
As $\delta<3$, it follows that
\[
\delta fq(q-1)
<3q^{5/2}<q^s<(q+1)^s.
\]

Since $\PSL_2(5)\cong\PSL_2(4)$, the case $q=5$
is covered by $q=4$.
We treat the cases $q=8$ and $q=16$ separately.
For the remaining prime powers $4\le q<64$,
the inequality $\delta fq(q^2-1)<k^{1+s}$
is verified in Appendix~\ref{app:rankone}.

\smallskip
\noindent
\textit{(ii) and (iii) $q=8$ or $q=16$.}
The preceding case analysis shows that the largest
soluble subgroups of $T$ are point stabilisers.
Since a point stabiliser of $T$ has order
$q(q-1)/d$, we obtain $\sigma(T)=56$ for $q=8$
and $\sigma(T)=240$ for $q=16$.

For soluble $B\le\Aut(T)$, we have
$[B:B\cap T]\le|\Out(T)|=f$.
Hence $|B|\le f\sigma(T)$.
Since $h=[T:B\cap T]\ge q+1$, we have
\[
\begin{aligned}
\delta|B|&\le168\delta
=9^{r_0}\le h^{r_0}
&&\text{if }q=8,\\
\delta|B|&\le960\delta
=17^{r_{16}}\le h^{r_{16}}
&&\text{if }q=16.
\end{aligned}
\]

Equality in either bound forces equality
throughout the preceding chain. Hence
$h=q+1$ and $|B|=f\sigma(T)$.
Thus $B\cap T$ is a point stabiliser of $T$.
Since $B$ normalises $B\cap T$, it fixes
the unique point $\alpha\in\mathbb P^1(\F_q)$
fixed by $B\cap T$.
Consequently, $B\le\Aut(T)_\alpha$.
As $|B|=fq(q-1)=|\Aut(T)_\alpha|$,
we have $B=\Aut(T)_\alpha$.
Conversely, each such point stabiliser
attains equality in the corresponding bound.
A direct calculation gives
$(960\delta)^5<17^{14}$ and $9^{14}<(168\delta)^5$,
so $r_{16}<14/5<r_0$.
\end{proof}

\begin{remark}\label{rem:rankone-exceptions}
Although \eqref{eq:outer-bound} fails for
$T\cong\PSL_2(8)$ and $\PSL_2(16)$,
Proposition~\ref{prop:rankone} gives
$
\delta|B|\le[T:B\cap T]^{r_0}
$
for every soluble $B\le\Aut(T)$,
since $r_{16}<r_0$.
\end{remark}

\begin{proposition}\label{prop:alt}
For every alternating simple group $T=A_n$, $n\ge5$, the inequality
\eqref{eq:outer-bound} holds.
\end{proposition}
\begin{proof}
Since $s=39/14$ and $\delta=24^{1/3}$,
taking $42$nd powers in \eqref{eq:outer-bound}
gives the equivalent inequality
\[
24^{14}|\Out(A_n)|^{42}\sigma(A_n)^{159}
<|A_n|^{117}.
\]

\smallskip
\noindent
\textit{(i) $5\le n\le9$.}
Since $A_5\cong\PSL_2(4)$ and
$A_6\cong\PSL_2(9)$,
the result for $n=5,6$ follows from
Proposition~\ref{prop:rankone}.
For $7\le n\le9$, Mann's Theorems~1 and~3
\cite{Mann} give $\sigma(A_n)\le m_n$, where
\begin{center}
\begin{minipage}{\linewidth}
\centering
\captionof{table}{Soluble subgroup bounds for small alternating groups}
\label{tab:alternating-small}
\[
\begin{array}{c|rrr}
n&7&8&9\\\hline
m_n&144&576&1296
\end{array}
\]
\end{minipage}
\end{center}
Since $|\Out(A_n)|=2$ and $|A_n|=n!/2$,
the bounds $\sigma(A_n)\le m_n$ show that
it suffices to verify
\[
24^{14}2^{42}m_n^{159}<(n!/2)^{117}.
\]
These inequalities hold for $n=7,8,9$,
proving \eqref{eq:outer-bound}.

\smallskip
\noindent
\textit{(ii) $n\ge10$.}
By Dixon's bound \cite[Theorem~3]{Dixon},
$\sigma(A_n)\le\delta^{n-1}$, and
$|\Out(A_n)|=2$.
Since $|A_n|=n!/2$ and
$\delta^{159}=24^{53}$,
it suffices to prove
\[
2^{42}24^{14+53(n-1)}<(n!/2)^{117}.
\]
For $n=10$, the left-hand side is $2^{42}24^{491}$.
Since $10!/2>4\cdot24^4$, the right-hand side exceeds
$2^{234}24^{468}$.
Thus the inequality follows from $24^{23}<2^{192}$.
When $n$ increases by one, the two sides
are multiplied by $24^{53}$ and
$(n+1)^{117}$, respectively.
Since $24<11^2$, we have
$24^{53}<11^{117}\le(n+1)^{117}$ for $n\ge10$.
Thus the inequality follows by induction.
\end{proof}

\begin{proposition}\label{prop:sporadic}
The inequality \eqref{eq:outer-bound} holds for
all $26$ sporadic simple groups and the Tits
group ${}^2F_4(2)'$.
\end{proposition}

\begin{proof}
Since $s-2(1+s)/3=11/42$ and
$\delta^{42}=24^{14}$, it suffices to prove
\[
\sigma(T)\le |T|^{2/3}
\qquad\text{and}\qquad
24^{14}|\Out(T)|^{42}<|T|^{11}.
\]
By \cite[Theorem~1]{Burness} for the $26$ sporadic
groups and \cite[Remark~1(d)]{Burness} for the
Tits group, the action of $T$ on $T/U$ has base
size at most three for every soluble $U\le T$.
Write $h=[T:U]$. Since the action is faithful,
a base of size at most three gives $|T|\le h^3$.
Hence $|U|=|T|/h\le |T|^{2/3}$, proving
$\sigma(T)\le |T|^{2/3}$.
Since $|\Out(T)|\le2$, $|T|\ge7920$ and
$24^{14}2^{42}<7920^{11}$, the second
inequality also holds.
\end{proof}

\section{Groups of Lie type: the Borel--Tits dichotomy}\label{sec:lie}

We now turn to simple groups of Lie type.
To establish \eqref{eq:outer-bound} for such a group $T$,
it suffices to prove
$\delta|\Out(T)|\,|U|^{1+s}<|T|^s$
for every soluble subgroup $U\le T$.
We distinguish whether $O_p(U)$ is trivial, where $p$ is the
defining characteristic.
If $O_p(U)\ne1$, the Borel--Tits theorem places $U$ in a
parabolic subgroup, and Proposition~\ref{prop:parabolic}
reduces the required inequality to the numerical
condition \eqref{eq:PIprime} below.
We also establish the bounds for the case $O_p(U)=1$
that are used in Sections~\ref{sec:classical}
and~\ref{sec:exceptional}.

Let $T$ be a simple group of Lie type in
characteristic $p$ over $\F_q$, where $q=p^f$.
We use this setup only for the classical groups
treated in Section~\ref{sec:classical} and for
$E_6(q)$, $E_7(q)$ and $E_8(q)$.
For classical $T$, take $X$ to be the
corresponding natural matrix group
$\SL_n(q)$, $\SU_n(q)$, $\Sp_n(q)$ or
$\Omega^\varepsilon_n(q)$.
For $T=E_6(q)$, $E_7(q)$ or $E_8(q)$,
let $\mathbf G$ be the simply connected simple algebraic
group of the same type over $\overline{\F}_p$,
and let $\sigma$ be the standard Frobenius endomorphism
of $\mathbf G$ over $\F_q$.
We then take $X$ to be the group of fixed points
of $\sigma$ in $\mathbf G$, so $X=\mathbf G^\sigma$.

In each case, $X/Z(X)\cong T$ and
$p\nmid |Z(X)|$.
Let $\pi:X\to T$ be the corresponding
quotient map. For $U\le T$, let
$\widehat U=\pi^{-1}(U)$.
Write $Z=Z(X)$ and $z=|Z|$. Then
\[
|\widehat U|=z|U|,
\qquad
O_p(\widehat U)\cong O_p(U).
\]

To prove the second assertion, write
$E=\pi^{-1}(O_p(U))$.
Since $E/Z$ is a $p$-group and $Z$ is
a central $p'$-subgroup, we have
$E=Z\times P$ for a Sylow $p$-subgroup
$P$ of $E$. As $P=O_p(E)$ and
$E\trianglelefteq\widehat U$, we have
$P\trianglelefteq\widehat U$.
Moreover, $O_p(\widehat U)\le E$, so
$O_p(\widehat U)=P$.
Since $\pi$ maps $P$ isomorphically onto $O_p(U)$,
the result follows.

\begin{theorem}[Borel--Tits
{\cite[Theorem~3.1.3(a)]{GLS3}}, \cite{BW}]
\label{thm:bt}
If $R$ is a nontrivial $p$-subgroup of $T$,
then there exists a proper parabolic subgroup
$P$ of $T$ such that $N_T(R)\le P$.
\end{theorem}

Let $U\le T$ be soluble with $O_p(U)\ne1$.
Since $U\le N_T(O_p(U))$, Theorem~\ref{thm:bt}
gives a proper parabolic subgroup of $T$ containing $U$,
and hence a maximal parabolic subgroup $P$ of $T$
such that $U\le P$.
Every parabolic subgroup of $T$ is the image
under $\pi$ of a parabolic subgroup of $X$.
Since every parabolic subgroup of $X$ contains
$Z$, the subgroup $\widehat P=\pi^{-1}(P)$
is a maximal parabolic subgroup of $X$.
Hence $\widehat U\le\widehat P$.

Write
\[
|O_p(\widehat P)|=q^u,\qquad
L=\widehat P/O_p(\widehat P),\qquad
N=[X:\widehat P].
\]
We call $L$ the \emph{Levi quotient} of $\widehat P$.
For the parabolic subgroups considered below,
$O_p(\widehat P)$ is the unipotent radical of $\widehat P$
(see Section~\ref{sec:classical} and the proof of
Lemma~\ref{lem:exceptional-split}),
so $L$ is isomorphic to a Levi complement of $\widehat P$.
An opposite unipotent radical has order $q^u$
and intersects $\widehat P$ trivially.
Its elements therefore represent distinct
cosets of $\widehat P$, so $N\ge q^u$.
The following proposition shows that,
if \eqref{eq:inner-bound} holds for every nonabelian
simple group of order less than $|T|$, then
every soluble $U\le T$ with $\widehat U\le\widehat P$
satisfies the required inequality, provided that
\begin{equation}\label{eq:PI}
\delta|\Out(T)|q^u c(L)<zN^s.
\end{equation}
Since $N\ge q^u$ and $\delta<3$,
condition \eqref{eq:PI} holds whenever
\begin{equation}\label{eq:PIprime}
3\bigl(|\Out(T)|/z\bigr)c(L)
\le q^{(s-1)u}.
\end{equation}

\begin{proposition}
\label{prop:parabolic}
Assume that \eqref{eq:inner-bound} holds
for every nonabelian simple group
of order less than $|T|$.
If \eqref{eq:PI} holds, then every soluble
$U\le T$ with $\widehat U\le\widehat P$
satisfies
$
\delta|\Out(T)|\,|U|^{1+s}<|T|^s.
$
\end{proposition}

\begin{proof}
Let $S_1,\ldots,S_m$ be the nonabelian
composition factors of $L$. Since $Z$ is abelian, these
are also composition factors of $\widehat P/Z\cong P$.
As $P<T$, we have $|S_i|<|T|$ for every $i$.
Since $\widehat U$ is soluble,
Lemma~\ref{lem:sigma} gives
$
z|U|=|\widehat U|
\le \sigma(\widehat P)
\le q^u c(L)\prod_{i=1}^m\sigma(S_i).
$
On the other hand, since
$|\widehat P|=q^u|L|$, we have
$
z|T|=N|\widehat P|
=Nq^u c(L)\prod_{i=1}^m|S_i|.
$
By the induction hypothesis
\eqref{eq:inner-bound},
$\sigma(S_i)^{1+s}<|S_i|^s$ for every $i$.
Consequently,
\[
\begin{aligned}
\delta|\Out(T)|\frac{|U|^{1+s}}{|T|^s}
\le
\frac{\delta|\Out(T)|q^u c(L)}{zN^s}
\prod_{i=1}^m
\frac{\sigma(S_i)^{1+s}}{|S_i|^s}
\le
\frac{\delta|\Out(T)|q^u c(L)}{zN^s}
<1,
\end{aligned}
\]
where the last inequality follows from
\eqref{eq:PI}.
\end{proof}

Thus, the case $O_p(U)\ne1$ reduces to
verifying \eqref{eq:PI} for the relevant
maximal parabolic subgroups.
We carry out these verifications in
Sections~\ref{sec:classical} and~\ref{sec:exceptional}.

We next consider the case $O_p(U)=1$.
The following lemma allows us to apply
bounds for completely reducible soluble
linear groups.

\begin{lemma}\label{lem:graded}
Let $A\le\GL_d(q)$, where $q=p^f$, and let
$V=\F_q^d$. 
Consider a composition series of
$V$ as an $\F_q A$-module:
$
0=V_0<V_1<\cdots<V_k=V.
$
Write $d_i=\dim_{\F_q}(V_i/V_{i-1})$, and let
$K$ be the kernel of the action of $A$ on
$
\bigoplus_{i=1}^k V_i/V_{i-1}.
$
Then
$K\le O_p(A)$ and $|K|\le q^{\sum_{i<j}d_id_j}$.
\end{lemma}

\begin{proof}
Choose a basis of $V$ adapted to the composition series
by successively extending a basis of each $V_i$.
Since every $V_i$ is $A$-invariant, the elements of $A$
are represented by upper block-triangular matrices,
whose diagonal blocks describe the induced actions
on $V_i/V_{i-1}$.

Since $K$ acts trivially on each $V_i/V_{i-1}$,
every element of $K$ has the form
\[
\begin{pmatrix}
I_{d_1} & * & \cdots & *\\
0 & I_{d_2} & \cdots & *\\
\vdots & & \ddots & \vdots\\
0 & \cdots & 0 & I_{d_k}
\end{pmatrix}.
\]
The upper block-unitriangular matrices of this form
constitute a group of order
$
q^{\sum_{i<j}d_id_j},
$
since the $(i,j)$-block has $d_id_j$ entries
over $\F_q$ for each $i<j$.
Thus $K$ is a $p$-group of order at most
$q^{\sum_{i<j}d_id_j}$.
Since $K\lhd A$, we have $K\le O_p(A)$.
\end{proof}

We shall use the following two results for
completely reducible soluble linear groups.
For a positive integer $m$ and a prime $r$,
let $m_r$ denote the largest power of $r$
dividing $m$.

\begin{theorem}[P\'alfy--Wolf {\cite{Palfy}, \cite[Theorem~3.1(a)]{Wolf}}]
\label{thm:wolf}
Let $V$ be a nonzero finite vector space and let
$G\le\GL(V)$ be soluble and completely reducible.
Then
$
|G|\le\delta^{-1}|V|^\beta,
$
where $\beta=\log(48\delta)/\log9<2.2441$.
\end{theorem}

\begin{theorem}[Dixon \cite{DixonNormal}]\label{thm:dixonlin}
Let $F$ be a perfect field, let $r$ be a prime,
and let $G\le\GL_n(F)$ be finite, soluble and
completely reducible. Then
$
|G:O_r(G)|_r\le r^{\lambda_r(n)},
$
where
\[
\lambda_r(n)=
\begin{cases}
\lfloor 4n/3\rfloor-1,
& r=2,\\[2pt]
\displaystyle\sum_{i\ge0}
\left\lfloor\frac{n}{(r-1)r^i}\right\rfloor,
& r\text{ is a Fermat prime},\\[8pt]
\displaystyle\sum_{i\ge1}
\left\lfloor\frac{n}{r^i}\right\rfloor,
& r\text{ is odd and not a Fermat prime}.
\end{cases}
\]
If $r=\operatorname{char}F$, then $O_r(G)=1$
and $|G|_r\le r^{\lambda_r(n)}$.
If $r>n+1$, then $\lambda_r(n)=0$
and $G$ has a normal Sylow $r$-subgroup.
\end{theorem}

For a prime $\ell$ and an integer $a$
coprime to $\ell$, let $\ord_\ell(a)$ denote
the multiplicative order of $a$ modulo $\ell$.
The following lemma applies these results
to soluble subgroups of classical groups
with trivial $p$-core.
We use these bounds in
Section~\ref{sec:classical} to establish
\eqref{eq:outer-bound}.

\begin{lemma}\label{lem:corefree}
Let $X$ be a classical group with natural module
$\F_Q^n$, where $Q=q^2$ for unitary groups
and $Q=q$ otherwise.
Let $U\le T$ be soluble with $O_p(U)=1$.
Then
\begin{enumerate}[label=\textup{(\alph*)}]
\item $|U|\le Q^{\beta n}/(\delta(Q-1))$;
\item $|U|_p\le p^{\lambda_p(n)}$;
\item If $\ell\ne p$ is a prime divisor of $|U|$
with $\ell>n+1$ and $j=\ord_\ell(Q)>n/2$,
then $|U|\le j(Q^j-1)|\GL_{n-j}(Q)|/(Q-1)$.

\end{enumerate}
\end{lemma}

\begin{proof}

Let $V=\F_Q^n$, and identify $\F_Q^\times$
with the subgroup of scalar matrices in $\GL(V)$.
Since $X\cap\F_Q^\times=Z\le\widehat U$,
the group $J=\widehat U\F_Q^\times$ is
soluble and has order $(Q-1)|U|$.
As $\F_Q^\times$ is central in $\GL(V)$,
we have $\widehat U\trianglelefteq J$.
Since $[J:\widehat U]$ is prime to $p$,
it follows that
$O_p(J)\le O_p(\widehat U)=1$.

By Lemma~\ref{lem:graded}, $J$ acts
faithfully on the direct sum of the
composition factors of $V$ as an
$\F_QJ$-module.
This module is completely reducible
and has dimension $n$.
Thus $J$ is isomorphic to a completely
reducible soluble subgroup
$J^*\le\GL_n(Q)$.

\smallskip
\noindent
\textit{(a).}
By Theorem~\ref{thm:wolf},
$(Q-1)|U|=|J|\le\delta^{-1}Q^{\beta n}$.

\smallskip
\noindent
\textit{(b).}
Applying Theorem~\ref{thm:dixonlin}
with $r=p$ gives
$|J|_p\le p^{\lambda_p(n)}$.
Since $p\nmid(Q-1)$, we have
$|U|_p=|J|_p$.

\smallskip
\noindent
\textit{(c).}
Since $\ell\mid|U|$ and
$J\le\GL_n(Q)$, we have
$\ell\mid|\GL_n(Q)|$, so $j\le n$.
As $j=\ord_\ell(Q)>n/2$,
the only factor $Q^k-1$ with
$1\le k\le n$ divisible by $\ell$
is $Q^j-1$.
Thus a Singer cyclic subgroup of $\GL_j(Q)$
contains a Sylow $\ell$-subgroup
of $\GL_n(Q)$, so the Sylow
$\ell$-subgroups of $J$ are cyclic.
Since $\ell>n+1$,
Theorem~\ref{thm:dixonlin} applied
to $J^*$ gives a normal Sylow
$\ell$-subgroup of $J$.
Since this Sylow $\ell$-subgroup is cyclic,
it has a unique subgroup $\langle x\rangle$
of order $\ell$.
As the Sylow subgroup is normal in $J$,
so is $\langle x\rangle$.

We now return to the action of $J$
on the natural module $V$.
Since $\ell\ne p$, the action of
$\langle x\rangle$ is completely
reducible, and every nontrivial
irreducible constituent has dimension $j$.
As $x\ne1$ and $2j>n$, there is exactly
one such constituent.
Write $[V,x]=(x-1)V$.
Thus $V=[V,x]\oplus C_V(x)$,
where $[V,x]$ is irreducible of dimension $j$.

Since $J$ normalises $\langle x\rangle$,
both $[V,x]$ and $C_V(x)$ are $J$-invariant.
The minimal polynomial of $x$ on $[V,x]$
is irreducible of degree $j$.
Hence the algebra $\F_Q[x]$ generated
by this restriction is a field
isomorphic to $\F_{Q^j}$.
The image of $J$ in $\GL([V,x])$
also normalises this field.
Identifying $[V,x]$ with $\F_{Q^j}$, we see that the normaliser of
$\F_Q[x]$ in $\GL([V,x])$ consists of the maps
$v\mapsto b v^{Q^i}$ with $b\in\F_{Q^j}^{\times}$
and $0\le i<j$, and has order $j(Q^j-1)$.
Finally, as $J$ acts faithfully on $V$,
we have
$(Q-1)|U|=|J|
\le j(Q^j-1)|\GL_{n-j}(Q)|$.

\end{proof}

We next establish bounds for the groups
$E_6(q)$, $E_7(q)$ and $E_8(q)$.
The case $E_6(2)$ requires the following
bound on the Weyl group.

\begin{lemma}\label{lem:oddW}
Every odd-order nilpotent subgroup of the Weyl group $W(E_6)$ has order at
most $81$.
\end{lemma}

\begin{proof}
Let $K$ be an odd-order nilpotent subgroup
of $W(E_6)$.
Since $|W(E_6)|=2^73^45$, if $5\nmid|K|$,
then $K$ is a $3$-group and $|K|\le81$.
Now suppose $5\mid|K|$.
Choose $x\in K$ of order $5$, and let
$A$ be a Sylow $3$-subgroup of $K$.
Since $K$ is nilpotent, $A$ centralises $x$.

Let $W$ be the six-dimensional rational
reflection module of $W(E_6)$.
Write $\Phi_5(t)=t^4+t^3+t^2+t+1$
for the fifth cyclotomic polynomial.
The characteristic polynomial of $x$
on $W$ is $(t-1)^2\Phi_5(t)$.
Thus $W=W_0\oplus W_1$, where $W_0$ is
the $2$-dimensional fixed space of $x$
and $W_1$ is a $4$-dimensional
$x$-invariant subspace.
The minimal polynomial of $x$ on $W_1$
is $\Phi_5(t)$.
Hence $W_1$ is one-dimensional over
$\mathbb Q(\zeta_5)$.
Since $A$ centralises $x$, it preserves
both $W_0$ and $W_1$.
Its action on $W_1$ is
$\mathbb Q(\zeta_5)$-linear, so its image
lies in $\mathbb Q(\zeta_5)^\times$.
The torsion subgroup of this group
has order $10$.
As $A$ is a $3$-group, it acts trivially
on $W_1$.
Since the reflection representation
is faithful, $A$ acts faithfully on $W_0$.

If $A\ne1$, choose $a\in Z(A)$ of order $3$.
As $A$ acts faithfully on $W_0$,
the element $a$ acts nontrivially on $W_0$.
Since $\dim_{\mathbb Q}W_0=2$,
its minimal polynomial is $t^2+t+1$.
Hence
$C_{\GL(W_0)}(a)\cong\mathbb Q(\zeta_3)^\times$.
The torsion subgroup of
$\mathbb Q(\zeta_3)^\times$ has order $6$.
Since $a\in Z(A)$, we have
$A\le C_{\GL(W_0)}(a)$, so $|A|\le3$.
This also holds when $A=1$.
Thus $|K|=5|A|\le15$.

\end{proof}

\begin{lemma}\label{lem:Ecorefree}
Let $T=E_6(q)$, $E_7(q)$ or $E_8(q)$, where $q=p^f$, and let
$U\le T$ be soluble with $O_p(U)=1$. Then
\[
 |U|<q^{52},\qquad |U|<q^{91},\qquad |U|<q^{104},
\]
respectively.
\end{lemma}

\begin{proof}

We may assume $U\ne1$, and let $F=F(U)$
be the Fitting subgroup of $U$.
Since $F$ is nilpotent, its Sylow
$p$-subgroup is characteristic in $F$
and hence normal in $U$.
As $O_p(U)=1$, it follows that
$F$ is a $p'$-group.

\medskip
\noindent
\textit{Step 1. Bounding $|F|$.}
Let $\mathbf G_{\mathrm{ad}}=\mathbf G/Z(\mathbf G)$
be the corresponding adjoint algebraic group.
The endomorphism $\sigma$ also induces
a Steinberg endomorphism of
$\mathbf G_{\mathrm{ad}}$.
We identify $T$ with a subgroup of
$\mathbf G_{\mathrm{ad}}^\sigma$.
Applying Vdovin~\cite[Lemma~2.3]{Vdovin} in
$\mathbf G_{\mathrm{ad}}^\sigma$, we see that $F$ normalises
a $\sigma$-stable maximal torus $\mathbf A$
of $\mathbf G_{\mathrm{ad}}$.
Write
$W=N_{\mathbf G_{\mathrm{ad}}}(\mathbf A)/\mathbf A$
for the Weyl group.
Since $F\le N_{\mathbf G_{\mathrm{ad}}}(\mathbf A)$,
the natural quotient map induces an embedding
$F/(F\cap\mathbf A)\hookrightarrow W$.
Thus $F/(F\cap\mathbf A)$ is a nilpotent
$p'$-subgroup of $W$.

Let $\ell=\dim\mathbf A$ be the rank of
$\mathbf G_{\mathrm{ad}}$.
Thus $\ell=6,7,8$ for the three groups,
respectively.
Let $D$ be the endomorphism induced by
$\sigma$ on the cocharacter lattice
of $\mathbf A$.
Since these groups are split over $\F_q$,
the endomorphism $D/q$ has finite order.
Hence $D^m=q^m I$ for some positive
integer $m$, so every eigenvalue of
$D$ has absolute value $q$.
The fixed-point formula for algebraic
tori now gives
$
|\mathbf A^\sigma|
=|\det(D-I)|\le(q+1)^\ell.
$
The kernel of the map $F\to W$ is
$F\cap\mathbf A\le\mathbf A^\sigma$.
Consequently,
$
|F|\le w(q+1)^\ell,
$
where $w$ is an upper bound for $|F/(F\cap\mathbf A)|$.
We may use $|W|_{p'}$, or $81$ for $T=E_6(2)$ by
Lemma~\ref{lem:oddW}.

Using the Weyl-group orders~\cite{Carter}
and Lemma~\ref{lem:oddW}, we obtain the
following choices of $w$.
The exponents $a$ are chosen to establish
$w(q+1)^\ell<q^a$.
\begin{center}
\begin{minipage}{\linewidth}
\centering
\captionof{table}{Weyl-group bounds for split exceptional groups}
\label{tab:weyl-bounds}
\[
\begin{array}{c|c|r|r}
T&q&w&a\\ \hline
E_6&2&81&16\\
E_6&3&640&16\\
E_6&\ge4&51840&16\\
E_7&2&2835&28\\
E_7&\ge3&2903040&28\\
E_8&2&42525&32\\
E_8&\ge3&696729600&32
\end{array}
\]
\end{minipage}
\end{center}
For example, $|W(E_6)|=2^7\cdot3^4\cdot5$,
so $w=640$ when $q=3$. In this case,
$640\cdot4^6<3^{16}$.
Since $a>\ell$, the ratio
$(q+1)^\ell/q^a$ decreases with $q$.
A direct calculation shows that $w(q+1)^\ell<q^a$
at the smallest value of $q$ in each row of
Table~\ref{tab:weyl-bounds}.
Hence $|F|<q^{16},q^{28},q^{32}$,
respectively.

\medskip
\noindent
\textit{Step 2. Bounding $|U|$.}
By Gasch\"utz's theorem~\cite[III.4.5]{Huppert},
$F/\Phi(U)=F(U/\Phi(U))$ is a direct product
of abelian minimal normal subgroups of $U/\Phi(U)$.
By \cite[III.4.2(b),(c)]{Huppert}, we have
$C_{U/\Phi(U)}(F/\Phi(U))\le F/\Phi(U)$
and $\Phi(U)<F$.
Hence $U/F$ acts faithfully on the nontrivial group
$F/\Phi(U)$.

Write $F/\Phi(U)=\prod_r A_r$ for its
primary summands.
For each nontrivial $A_r$, let $G_r$
be the image of $U/F$ on $A_r$.
Then $G_r$ is soluble and completely
reducible, so Theorem~\ref{thm:wolf} gives
$|G_r|\le\delta^{-1}|A_r|^\beta$.
Since $U/F$ embeds in the product of
these groups, and $\delta>1$, we have
$
|U/F|\le\delta^{-1}|F/\Phi(U)|^\beta
\le\delta^{-1}|F|^\beta.
$
As $\beta<9/4$, it follows that
$
|U|\le\delta^{-1}|F|^{1+\beta}
<|F|^{13/4}.
$
The three bounds on $|F|$ now give
the required bounds on $|U|$.
\end{proof}

\section{Classical groups}\label{sec:classical}

In Section~\ref{sec:lie}, we reduced the case
$O_p(U)\ne1$ to bounds for parabolic subgroups
and obtained estimates for the case $O_p(U)=1$.
We now use these results to establish
\eqref{eq:outer-bound} for classical simple groups,
apart from $\PSL_3(3)$, which is treated separately.

The finite simple classical groups arise from
four families: linear, unitary, symplectic
and orthogonal.
These are represented by $\PSL_n(q)$,
$\PSU_n(q)$, $\PSp_n(q)$, and
$\Omega_n(q)$ or $\PO_n^\pm(q)$,
respectively, with suitable restrictions
on the parameters.
Since the groups $\PSL_2(q)$ have already
been treated in Section~\ref{sec:small},
we use the following standard low-dimensional
isomorphisms \cite[Section~2.9]{KL} to avoid
repetition:
\begin{center}
\begin{minipage}{\linewidth}
\centering
\captionof{table}{Low-dimensional classical-group isomorphisms}
\label{tab:classical-isomorphisms}
\begin{tabular}{ll}
\toprule
Group & Isomorphic group \\ \midrule
$\PSU_2(q),\ \PSp_2(q)$
  & $\PSL_2(q)$ \\
$\Omega_3(q)$, $q$ odd
  & $\PSL_2(q)$ \\
$\PO_4^+(q)$
  & $\PSL_2(q)\times\PSL_2(q)$ \\
$\PO_4^-(q)$
  & $\PSL_2(q^2)$ \\
$\Omega_5(q)$, $q$ odd
  & $\PSp_4(q)$ \\
$\PO_6^+(q)$
  & $\PSL_4(q)$ \\
$\PO_6^-(q)$
  & $\PSU_4(q)$ \\
$\Omega_{2m+1}(q)$, $q$ even
  & $\PSp_{2m}(q)$ \\
\bottomrule
\end{tabular}
\end{minipage}
\end{center}

We also exclude $\PSL_3(2)\cong\PSL_2(7)$
and $\PSL_4(2)\cong A_8$, which have
already been treated, together with
the soluble $\PSU_3(2)$ and the nonsimple
$\PSp_4(2)\cong S_6$.

Thus, apart from the four small cases excluded above,
the remaining classical groups are as follows:
\begin{table}[htbp]
\centering
\caption{Remaining classical groups}
\label{tab:classical-cases}
\small
\begin{tabular}{@{}lll@{}}
\toprule
type & group & parameters\\
\midrule
Linear
  & $\PSL_n(q)$ & $n\ge3$\\
Unitary
  & $\PSU_n(q)$ & $n\ge3$\\
Symplectic
  & $\PSp_n(q)$ & $n\ge4$ even\\
Orthogonal
  & $\Omega_n(q)$ & $n\ge7$ odd, $q$ odd\\
  & $\PO_n^\pm(q)$ & $n\ge8$ even\\
\bottomrule
\end{tabular}
\end{table}

We write $U_n(q)=\PSU_n(q)$ and
$S_n(q)=\PSp_n(q)$.
In the tables below, $L_n$, $U_n$, $S_n$,
$O_n$ and $O_n^\pm$ are used as labels
for the linear, unitary, symplectic,
odd-dimensional orthogonal and even-dimensional
orthogonal types, respectively.
We use the following standard data
\cite[Table~5.1.A]{KL}.
Here $q=p^f$, where $p$ is prime, $z=|Z(X)|$, $|T|_p=q^\nu$,
and $D$ is the dimension of the
corresponding simple algebraic group.
For even-dimensional orthogonal groups,
$\varepsilon\in\{+1,-1\}$ denotes the
plus or minus type.

\begin{center}
\begin{minipage}{\linewidth}
\centering
\captionof{table}{Parameters of classical groups}
\label{tab:classical-parameters}
\small\setlength{\tabcolsep}{2pt}
\begin{tabular}{lcccc}\toprule
type & $z$ & $|\Out(T)|$ & $\nu$  & $D$\\\midrule
$L_n$ & $\gcd(n,q-1)$ & $2f\gcd(n,q-1)$ & $n(n-1)/2$ & $n^2-1$\\
$U_n$ & $\gcd(n,q+1)$ & $2f\gcd(n,q+1)$ & $n(n-1)/2$ & $n^2-1$\\
$S_n$ & $\gcd(2,q-1)$ & $f\gcd(2,q-1)$; $2f$ if $n=4$, $p=2$ & $n^2/4$ & $n(n+1)/2$\\
$O_n$, $n$ odd & $1$ & $2f$ & $(n-1)^2/4$ & $n(n-1)/2$\\
$O^\varepsilon_n$, $n=2m$ & $\gcd(4,q^m-\varepsilon)/\gcd(2,q-1)$ & $2f\gcd(4,q^m-\varepsilon)$ & $m(m-1)$ & $n(n-1)/2$\\
\bottomrule
\end{tabular}
\end{minipage}
\end{center}
For $O_8^+(q)$, the outer automorphism group
has order $6f\gcd(4,q^4-1)$ instead.
In all cases, $|\Out(T)|\le12fz$.
The standard order formulas also give
$q^D/(4z)\le|T|<q^D$.

\begin{proposition}\label{prop:linear}
Let $T=\PSL_n(q)$ with $n\ge3$. Then \eqref{eq:outer-bound} holds
unless $(n,q)=(3,3)$.
\end{proposition}

\begin{proof}
\noindent
\textit{(1) Determining $\sigma(T)$.}
Let $H$ be a maximum-order soluble subgroup of $\SL_n(q)$.
Since $n\ge3$, Mann's Theorem~8 \cite{Mann} gives a
maximum-order soluble subgroup $M$ of $\GL_n(q)$
with $H=M\cap\SL_n(q)$.
By Mann's Theorems~5 and~6, in a suitable basis
$M$ is triangular when $q=4$ or $q\ge7$.
For $q=2,3,5$, it is block triangular, with $\lfloor n/2\rfloor$
diagonal blocks of dimension two and
one block of dimension one when $n$ is odd.
The two-dimensional blocks have orders
$6$, $48$ and $96$ for $q=2,3,5$,
respectively.

Write
\[
c_q=
\begin{cases}
3,&q=2,\\
4,&q=3,\\
6/5,&q=5,\\
1,&q=4\text{ or }q\ge7.
\end{cases}
\]
Counting the blocks gives
$|M|=q^{n(n-1)/2}(q-1)^n
c_q^{\lfloor n/2\rfloor}$.

To pass from $\GL_n(q)$ to $\SL_n(q)$,
we first show that $\det(M)=\F_q^\times$,
where
$\det:\GL_n(q)\to\F_q^\times$
is the determinant homomorphism.
This is clear when $M$ is triangular.
For $q=2$, the group $\F_q^\times$ is
trivial, while for $q=3$ the
two-dimensional blocks are $\GL_2(3)$.
For $q=5$, let $B\le\GL_2(5)$ be a
diagonal block of order $96$, and write
$e=|\det(B)|$.
Since $e\mid4$ and
$96/e$ divides $|\SL_2(5)|=120$,
we have $e=4$.
Thus $\det(M)=\F_q^\times$ in every case.

It follows that $|H|=|M|/(q-1)$.
Write $Z=Z(\SL_n(q))$, so
$|Z|=\gcd(n,q-1)=d$.
Since $HZ$ is soluble and $H$ has
maximum order, we have $Z\le H$.

Let $\pi:\SL_n(q)\to T$ be the natural
quotient map with kernel $Z$.
Since $Z$ is abelian, $\pi^{-1}(U)$
is soluble for every soluble subgroup
$U\le T$.
As $H$ has maximum order, we obtain
$d|U|=|\pi^{-1}(U)|\le|H|$.
On the other hand, $\pi(H)$ is soluble
and has order $|H|/d$.
Hence
\begin{equation}\label{eq:mann-sigma}
\sigma(T)=
\frac{q^{n(n-1)/2}(q-1)^{n-1}}{d}
c_q^{\lfloor n/2\rfloor}.
\end{equation}
We now prove \eqref{eq:outer-bound}.
Since
\[
|T|=\frac{q^{n^2-1}}{d}
\prod_{i=2}^n(1-q^{-i}),
\qquad |\Out(T)|=2df,
\]
equation \eqref{eq:mann-sigma} shows that
\eqref{eq:outer-bound} is equivalent to
\[
\frac{2\delta f\,c_q^{53\lfloor n/2\rfloor/14}
(1-q^{-1})^{53(n-1)/14}}
{q^{(n-1)(25n-28)/28}
\left(\prod_{i=2}^n(1-q^{-i})\right)^{39/14}}
<1.
\]
Write
\[
P_n(q)=\prod_{i=2}^n(1-q^{-i}),
\qquad
a_n=\frac{(n-1)(25n-28)}{28},
\]
and define
\[
R_{n,q}=
\frac{2\delta f\,c_q^{53\lfloor n/2\rfloor/14}
(1-q^{-1})^{53(n-1)/14}}
{q^{a_n}P_n(q)^{39/14}}.
\]
It remains to prove $R_{n,q}<1$.
We consider $q\ge4$ and $q=2,3$
separately.

\medskip
\noindent
\textit{(2) The case $q\ge4$.}
We have $c_q\le6/5$ and
$P_n(q)>1-1/(q(q-1))\ge11/12$.
Also, $a_n\ge n$ for $n\ge3$, since
$28(a_n-n)=25n^2-81n+28>0$.

As $53/14<4$ and $39/14<3$, we have
\[
(6/5)^{53/14}<(6/5)^4<(3/2)^2,
\qquad
(12/11)^{39/14}<(12/11)^3<4/3.
\]
Since $\delta<3$ and $f\le q/2$,
it follows that
\[
R_{n,q}
<8f q^{-a_n}(3/2)^n
\le4q\left(\frac{3}{2q}\right)^n
\le\frac{27}{2q^2}
\le\frac{27}{32}<1.
\]

\medskip
\noindent
\textit{(3) The cases $q=2,3$.}
We show that $R_{n,q}$ decreases with $n$.
Write
$\epsilon_n=\lfloor(n+1)/2\rfloor-\lfloor n/2\rfloor$.
Since $n\ge3$, we have
$(25n-14)/14>4$ and
$1-q^{-(n+1)}\ge1-q^{-4}$.
Moreover, $s<3$ and
$c_q(1-q^{-1})>1$.
Thus
\[
\begin{aligned}
\frac{R_{n+1,q}}{R_{n,q}}
=
\frac{c_q^{53\epsilon_n/14}
(1-q^{-1})^{53/14}}
{q^{(25n-14)/14}
(1-q^{-(n+1)})^{39/14}}
<
\frac{[c_q(1-q^{-1})]^4}
{q^4(1-q^{-4})^3}
=
\begin{cases}
48/125,&q=2,\\
81/125,&q=3.
\end{cases}
\end{aligned}
\]
Hence $R_{n,q}$ is strictly decreasing
with $n$ for both values of $q$.

Finally, direct calculation gives
\[
R_{3,2}^{42}
=\frac{2^{210}3^{56}}{7^{117}}<1,
\qquad
R_{3,3}^{42}
=\frac{2^{252}3^{140}}{13^{117}}>1,
\]
and
\[
R_{4,3}^{42}
=\frac{2^{261}3^{266}}{65^{117}}<1.
\]
Therefore $R_{n,q}<1$ for all $n\ge3$
and $q=2,3$, except when $(n,q)=(3,3)$.
Together with (2), this proves
the proposition.
\end{proof}

\begin{proposition}\label{prop:L33}
Let $T=\PSL_3(3)$. Then
$
\delta|B|<[T:B\cap T]^s
$
for every soluble subgroup $B\le\Aut(T)$.
\end{proposition}

\begin{proof}
Here $|T|=5616$ and $\sigma(T)=432$ by
\eqref{eq:mann-sigma}.
The maximal subgroups of $\Aut(T)$ other than
$T$ have orders $216,96,78,48$ (see \cite{ATLAS}).
Hence every soluble subgroup $B\le\Aut(T)$
has order at most $432$.
Write $h=[T:B\cap T]$. Since $|B\cap T|\le432$,
we have $h\ge13$. The integer comparison
$
24^{14}\cdot432^{42}<13^{117}
$
gives
$
\delta|B|\le\delta\cdot432<13^s\le h^s
$.
\end{proof}

\begin{remark}\label{rem:L33-exception}
Although \eqref{eq:outer-bound} fails for
$T\cong\PSL_3(3)$,
Proposition~\ref{prop:L33} gives
$
\delta|B|<[T:B\cap T]^{r_0}
$
for every soluble $B\le\Aut(T)$,
since $s<r_0$.
\end{remark}

\begin{lemma}\label{lem:classindex}
Let $T$ be one of the following groups:
\[
U_3(q)\ (q\ge3),\quad
U_5(q)\ (q\ge2),\quad
S_4(q)\ (q\ge3),\quad
U_4(2).
\]
For every proper subgroup $U<T$,
write $h=[T:U]$. Then
$
\delta|\Out(T)|\,|T|<h^{15/4}.
$
\end{lemma}

\begin{proof}
Write $t=|T|$ and $o=|\Out(T)|$.
Let $\mu(T)$ denote the least index of
a proper subgroup of $T$.
Since $h\ge\mu(T)$, it suffices to prove
$
\delta ot<\mu(T)^{15/4}.
$
Since $T$ is nonabelian simple, every proper
subgroup is core-free.
Hence $\mu(T)$ equals the minimum degree
of a faithful permutation representation
of $T$, as tabulated in \cite[Table~4]{GMPS}.

\medskip
\noindent
\textit{(1) The groups $U_3(q)$.}
For $q\ne5$, we have $\mu(T)=q^3+1$ and
$
ot=2fq^3(q^3+1)(q^2-1).
$
Since $\delta<3$ and $f\le q/2$, we have
\[
\frac{\delta ot}{\mu(T)^{15/4}}
<\frac{6fq^5}{q^{33/4}}
\le3q^{-9/4}<1.
\]
\medskip
\noindent
\textit{(2) The groups $U_5(q)$.}
We have
\[
\mu(T)=\frac{(q^5+1)(q^4-1)}{q^2-1}
=(q^5+1)(q^2+1)>q^7.
\]
Since
$(1-q^{-2})(1+q^{-3})<1$ and
$(1-q^{-4})(1+q^{-5})<1$,
the order formula gives
\[
ot=2fq^{10}(q^2-1)(q^3+1)
(q^4-1)(q^5+1)<2fq^{24}.
\]
For $q\ge3$, we have $f\le q/2$
and $3<q^{5/4}$.
Thus
$
\delta ot<6fq^{24}
\le3q^{25}
<q^{105/4}
<\mu(T)^{15/4}.
$

\medskip
\noindent
\textit{(3) The groups $S_4(q)$.}
For $q>3$, 
we have $\mu(T)=(q^4-1)/(q-1)>q^3$.
If $q>3$ is odd, then $ot<fq^{10}$
and $3f\le q$. Hence
$
\delta ot<3fq^{10}
\le q^{11}
<q^{45/4}
<\mu(T)^{15/4}.
$

For $q=2^f$ with $f\ge4$, we have
$ot<2fq^{10}$.
Since $6\cdot4<2^5$ and
$(f+1)/f\le5/4<2^{5/4}$ for $f\ge4$,
induction gives $6f<q^{5/4}$.
Thus
$
\delta ot<6fq^{10}
<q^{45/4}
<\mu(T)^{15/4}.
$

\medskip
\noindent
\textit{(4) The remaining cases.}
Thus only $q=3,4,8$ remain for $S_4(q)$.
It remains to check $U_3(5)$, $U_5(2)$,
$S_4(3)$, $S_4(4)$ and $S_4(8)$.
The group $U_4(2)$ is also covered,
since $U_4(2)\cong S_4(3)$.
The values of $t$ and $o$ follow from
the standard formulas above, while
those of $\mu(T)$ are given in
\cite[Table~4]{GMPS}:
\begin{center}
\begin{minipage}{\linewidth}
\centering
\captionof{table}{Indices of small classical groups}
\label{tab:classical-small-indices}
\[
\begin{array}{c|r|r|r}
T&t&o&\mu(T)\\ \hline
U_3(5)&126000&6&50\\
U_5(2)&13685760&2&165\\
S_4(3)\cong U_4(2)&25920&2&27\\
S_4(4)&979200&4&85\\
S_4(8)&1056706560&6&585
\end{array}
\]
\end{minipage}
\end{center}
Since $\delta<3$, it suffices to verify
$
(3ot)^4<\mu(T)^{15}.
$
This holds for every row by direct
integer comparison.
\end{proof}

Combining Table~\ref{tab:classical-cases} with
Propositions~\ref{prop:linear} and~\ref{prop:L33}
and Lemma~\ref{lem:classindex}, we obtain the following summary:
\begin{center}
\begin{minipage}{\linewidth}
\centering
\captionof{table}{Completed and remaining classical cases}
\label{tab:classical-summary}
\small
\begin{tabular}{@{}lp{0.34\linewidth}p{0.40\linewidth}@{}}
\toprule
group & already treated & remaining cases\\
\midrule
$\PSL_n(q)$
  & All cases (Propositions~\ref{prop:linear}
    and~\ref{prop:L33})
  & ---\\
\addlinespace
$\PSU_n(q)$
  & Cases in Lemma~\ref{lem:classindex}
  & $U_4(q)$ with $q\ge3$;
    $U_n(q)$ with $n\ge6$\\
\addlinespace
$\PSp_n(q)$
  & Cases in Lemma~\ref{lem:classindex}
& $n\ge6$ even\\
\addlinespace
$\Omega_n(q)$
  & ---
  & $n\ge7$ odd, $q$ odd\\
\addlinespace
$\PO_n^\varepsilon(q)$
  & ---
  & $n\ge8$ even\\
\bottomrule
\end{tabular}
\end{minipage}
\end{center}
We now consider the remaining unitary, symplectic
and orthogonal groups, beginning with soluble
subgroups $U\le T$ satisfying $O_p(U)\ne1$.

Recall from Section~\ref{sec:lie} that
$X/Z(X)\cong T$, where $z=|Z(X)|$.
By the Borel--Tits theorem, every soluble
$U\le T$ with $O_p(U)\ne1$ is contained in
the image of a maximal parabolic subgroup
$\widehat P$ of $X$.
Write $|O_p(\widehat P)|=q^u$, and let
$L=\widehat P/O_p(\widehat P)$ be the
Levi quotient.
Since $[X:\widehat P]\ge q^u$ and $\delta<3$,
Proposition~\ref{prop:parabolic} shows that,
assuming \eqref{eq:inner-bound} for smaller
nonabelian simple groups, it suffices to verify
\eqref{eq:PIprime}, namely
\[
3\frac{|\Out(T)|}{z}c(L)\le q^{(s-1)u}.
\]
We therefore need the value of $u$ and an upper
bound for $c(L)$ for each maximal parabolic subgroup.
For the remaining classical groups, these subgroups
are stabilisers of totally singular $k$-spaces in
the unitary and orthogonal cases, and of totally
isotropic $k$-spaces in the symplectic case.

Recall that for any finite group $Y$, $c(Y)$
denotes the product of the orders of its abelian
composition factors.
To bound $c(L)$, we use the following notation:
$c_{\SL}(r,q)=c(\SL_r(q))$,
$c_{\SU}(r,q)=c(\SU_r(q))$ and
$c_{\Sp}(r,q)=c(\Sp_r(q))$.
For the full orthogonal groups, write
$c_{\mathrm O}(r,q)=c(\mathrm O_r(q))$ when $r$ is odd,
and $c_{\mathrm O}(r,q)=c(\mathrm O_r^\varepsilon(q))$
when $r>0$ is even, with $c_{\mathrm O}(0,q)=1$.
The sign $\varepsilon$ is determined by the
orthogonal type. 
The following values follow from the standard
orders and composition factors of the classical
groups (see \cite[Chapter~2]{KL}), together with
the low-dimensional isomorphisms listed above:
\begin{center}
\begin{minipage}{\linewidth}
\centering
\captionof{table}{Abelian composition-factor products for classical groups}
\label{tab:classical-c}
\small
\begin{tabular}{lll}
\toprule
function & case & value\\
\midrule
$c_{\SL}(r,q)$
  & $r\le1$ & $1$\\
  & $(r,q)=(2,2),(2,3)$ & $6,24$\\
  & otherwise & $\gcd(r,q-1)$\\
\addlinespace
$c_{\SU}(r,q)$
  & $r\le1$ & $1$\\
  & $(r,q)=(2,2),(2,3),(3,2)$ & $6,24,216$\\
  & otherwise & $\gcd(r,q+1)$\\
\addlinespace
$c_{\Sp}(r,q)$
  & $r=0$ & $1$\\
  & $(r,q)=(2,2),(2,3),(4,2)$ & $6,24,2$\\
  & otherwise & $\gcd(2,q-1)$\\
\addlinespace
$c_{\mathrm O}(r,q)$
  & $r=0$ & $1$\\
  & $r=1$ & $\gcd(2,q-1)$\\
  & $r=2$ & $2(q-\varepsilon)$\\
  & $(r,q)=(3,3)$ & $48$\\
  & $r\ge3$ odd, otherwise & $4$\\
  & $(r,q,\varepsilon)=(4,2,+1),(4,3,+1)$
  & $72,1152$\\
  & $r\ge4$ even, otherwise
  & $2\gcd(4,q^{r/2}-\varepsilon)$\\
\bottomrule
\end{tabular}
\end{minipage}
\end{center}
For the symplectic groups, $r$ is even;
in the odd-dimensional orthogonal case, $q$ is odd.

The standard Levi decompositions, together with
the determinant condition in the unitary case,
give the following values of $u$ and upper bounds
for $c(L)$, where $m=n-2k$:
\begin{table}[htbp]
\centering
\caption{Unipotent exponents and bounds for Levi quotients}
\label{tab:classical-levi}
\small
\begin{tabular}{lccc}
\toprule
type & $k$ & $u$ & upper bound for $c(L)$\\
\midrule
$U_n$
  & $1\le k\le n/2$
  & $k(k+2m)$
  & $(q^2-1)c_{\SL}(k,q^2)c_{\SU}(m,q)$\\
$S_n$
  & $1\le k\le n/2$
  & $km+k(k+1)/2$
  & $(q-1)c_{\SL}(k,q)c_{\Sp}(m,q)$\\
$O^{(\varepsilon)}_n$
  & \shortstack[l]{$1\le k\le$ Witt index\\
    ($k\ne n/2-1$ for $O_n^+$)}
  & $km+k(k-1)/2$
  & $(q-1)c_{\SL}(k,q)c_{\mathrm O}(m,q)$\\
\bottomrule
\end{tabular}
\end{table}
Let $Q$ be the unipotent radical of $\widehat P$.
In each row the usual Levi complement is normal in the
corresponding direct product of general linear and full isometry groups.
The factors of these products have trivial $p$-core.
Indeed, the general linear, general unitary and symplectic groups
act irreducibly on their natural modules, as do the full orthogonal
groups in dimension at least three. The factors of dimension zero
or one have order prime to $p$.
For an orthogonal factor of dimension two, the restriction on $k$
excludes the plus type. Its order is $2(q+1)$, and when $p=2$
it is dihedral with a cyclic subgroup of odd order $q+1>1$,
so its $2$-core is trivial.
Hence the Levi complement has trivial $p$-core, and
$O_p(\widehat P)=Q$. In particular, $L$ is isomorphic to the
Levi complement, and the values of $u$ are those of the
unipotent radical.
In the unitary row, the bound for $c(L)$ is exact
when $m>0$. When $m=0$, we have
$c(L)=c(\GL_k(q^2))/(q+1)$, so the tabulated
bound exceeds $c(L)$ by a factor of $q+1$.
In the orthogonal row, the Levi complement is the intersection of
$X$ with the block diagonal subgroup
$\GL_k(q)\times\mathrm O_m^{(\varepsilon)}(q)$ of the full isometry group,
where the second factor has the indicated type and is trivial for $m=0$.
Since $X$ is normal in the full isometry group, this intersection
is normal in the direct product. Hence
$c(L)\le c(\GL_k(q))c_{\mathrm O}(m,q)$.
In type $O_n^+$ with even $n$, the maximal totally
singular $n/2$-spaces form two $X$-orbits, whose
stabilisers have the same index in $X$ and the
same values of $u$ and $c(L)$.

The following lemma deals with the case $O_p(U)\ne1$.

\begin{lemma}\label{lem:classPI}
Let $T$ be one of the following groups:
\[
\begin{gathered}
U_4(q)\ (q\ge3),\quad
U_n(q)\ (n\ge6),\quad
S_n(q)\ (n\ge6\text{ even}),\\
\Omega_n(q)\ (n\ge7\text{ odd},\ q\text{ odd}),\quad
\PO_n^\varepsilon(q)\ (n\ge8\text{ even}).
\end{gathered}
\]
Assume that \eqref{eq:inner-bound} holds for every
nonabelian simple group of order less than $|T|$.
Let $U\le T$ be soluble with $O_p(U)\ne1$,
and write $h=[T:U]$. Then
$
\delta|\Out(T)|\,|T|<h^{1+s}.
$
\end{lemma}

\begin{proof}
By Theorem~\ref{thm:bt}, the full preimage
$\widehat U$ of $U$ in $X$ lies in a maximal
parabolic subgroup $\widehat P$ of $X$.
By Proposition~\ref{prop:parabolic},
it suffices to verify \eqref{eq:PIprime}
for every maximal parabolic subgroup of $X$.
Let $o=|\Out(T)|$. Since $s-1=25/14$,
\eqref{eq:PIprime} takes the form
$
3(o/z)c(L)\le q^{25u/14}.
$

\medskip
\noindent
\textit{(1) The case $q\ge4$.}
Table~\ref{tab:classical-c} gives
\[
\begin{gathered}
c_{\SL}(k,q)\le q-1,\qquad
c_{\SL}(k,q^2)\le q^2-1,\\
c_{\SU}(m,q)\le q+1,\qquad
c_{\Sp}(m,q)\le2,\qquad
c_{\mathrm O}(m,q)\le2(q+1).
\end{gathered}
\]
Together with Table~\ref{tab:classical-levi},
these bounds give $c(L)\le(q+1)^5$.
Since $o/z\le12f$ and $f\le q/2$, we have
$
3(o/z)c(L)\le18q(q+1)^5<55q^6.
$
If $u\ge5$, then $25u/14-6\ge41/14$.
Since $55^{14}<4^{41}$ and $q\ge4$, it follows that
$55q^6<q^{25u/14}$.

By Table~\ref{tab:classical-levi},
we have $u\ge5$ except for $U_4(q)$ with
$k=2$, where $m=0$ and $u=4$.
In this case, the same table gives
$c(L)\le2(q^2-1)$, while
$o/z=2f\le q$. Hence
$
3(o/z)c(L)\le6q^3<q^5<q^{100/14}
$.

\medskip
\noindent
\textit{(2) The cases $q=2,3$.}
For $q=2$, Tables~\ref{tab:classical-c}
and~\ref{tab:classical-levi} give
$c(L)\le1944$, while
Table~\ref{tab:classical-parameters} gives
$o/z\le6$.
For $q=3$, the corresponding bounds are
$c(L)\le55296$ and $o/z\le12$.
Here $1944=3\cdot3\cdot216$ and
$55296=2\cdot24\cdot1152$, using the
unitary and orthogonal rows, respectively.
Hence
\[
3(o/z)c(L)\le
\begin{cases}
3\cdot6\cdot1944=34992,&q=2,\\
3\cdot12\cdot55296=1990656,&q=3.
\end{cases}
\]
Since
\[
34992<2^{16}<2^{225/14},
\qquad
1990656<3^{14}<3^{200/14},
\]
\eqref{eq:PIprime} holds for $u\ge9$ when $q=2$
and for $u\ge8$ when $q=3$.

It remains to consider $u\le8$ when $q=2$
and $u\le7$ when $q=3$.
For fixed $n$, each formula for $u$ in
Table~\ref{tab:classical-levi} is a concave quadratic in $k$.
Evaluating it at the smallest and largest permitted values
of $k$ gives $u\ge n-2$ in every case.
Hence $n\le10$ when $q=2$ and
$n\le9$ when $q=3$.

\medskip
\noindent
\textit{(3) The remaining cases.}
For $U_n(q)$ with $n\ge6$,
Table~\ref{tab:classical-levi} gives $u\ge9$,
so these groups are already covered.
The enumeration in Appendix~\ref{app:classical-parabolic}
gives the remaining cases listed below,
grouped by $(q,u)$.
For each case, the tables above give $o/z$
and an upper bound for $c(L)$.
The last column gives the maximum $M$ of the
resulting bounds for $3(o/z)c(L)$ in each row.
\begin{center}
\begin{minipage}{\linewidth}
\centering
\captionof{table}{Remaining classical parabolic cases}
\label{tab:classical-parabolic-cases}
\[
\begin{array}{c|c|l|r}
q & u & (\text{type},k) & M\\
\hline
2&5&(S_6,1)&6\\
2&6&(S_6,3),(O_8^+,1),(O_8^+,4),(O_8^-,1)&36\\
2&7&(S_6,2),(S_8,1)&108\\
2&8&(O_{10}^{\pm},1)&12\\
\hline
3&4&(U_4,2)&96\\
3&5&(U_4,1),(S_6,1),(O_7,1)&1152\\
3&6&(S_6,3),(O_7,3),(O_8^+,1),
       (O_8^+,4),(O_8^-,1)&288\\
3&7&(S_6,2),(S_8,1),(O_7,2),(O_9,1)&13824
\end{array}
\]
\end{minipage}
\end{center}
Here $O_{10}^{\pm}$ includes both types.
For example, when $q=3$, the case
$O_7(3)$ with $k=2$ has $u=7$.
The tables give $o/z=2$ and
$c(L)\le2\cdot24\cdot48$, so
$
3(o/z)c(L)\le
3\cdot2\cdot(2\cdot24\cdot48)=13824,
$
which is the maximum for $(q,u)=(3,7)$.
Every value $M$ in the table satisfies
$M<q^{u+2}$. Since $u\ge4$, we have
$u+2<25u/14$. 
Hence
$
3(o/z)c(L)\le M<q^{u+2}<q^{25u/14}.
$
Thus \eqref{eq:PIprime} holds for every
maximal parabolic subgroup of $X$.
By Proposition~\ref{prop:parabolic},
$
\delta o|U|^{1+s}<|T|^s.
$
Since $|T|=h|U|$, this is equivalent to
$\delta o|T|<h^{1+s}$.
\end{proof}

\begin{lemma}\label{lem:classCF}
Let $T$ be as in Lemma~\ref{lem:classPI}, and let
$U\le T$ be soluble with $O_p(U)=1$.
Write $h=[T:U]$. Then
$
\delta|\Out(T)|\,|T|<h^{15/4}.
$
\end{lemma}

\begin{proof}
Write $t=|T|$, $o=|\Out(T)|$, and $Q=q^a$, where
$a=2$ in unitary type and $a=1$ otherwise.
We aim to prove $\delta ot<h^{15/4}$.
We first establish bounds from the order formulas.
We then use Lemma~\ref{lem:corefree}(a) to handle
the higher-dimensional groups and $O_7(q)$.
The remaining groups, $U_4(q)$ and $U_6(q)$,
are treated using Lemma~\ref{lem:corefree}(b),(c).

\medskip
\noindent
\textit{Step 1. Preliminary bounds.}
We first establish
\[
\delta zt>q^D,\qquad
ot<6fq^D,\qquad
z\le Q-1,
\]
using the standard order formulas for classical
groups and the values of $z$ and $o$ in
Table~\ref{tab:classical-parameters}.

We use $5/2<\delta<3$.
In unitary type, since $zt=|\SU_n(q)|$ and
$D=n^2-1$, the standard order formula gives
\[
\frac{zt}{q^D}
=\prod_{i=2}^n(1-(-1)^iq^{-i})
\ge\prod_{i=2}^n(1-q^{-i})
\ge1-\sum_{i=2}^nq^{-i}>\frac12.
\]
In symplectic type, writing $n=2m$, we have
\[
\frac{zt}{q^D}
=\prod_{i=1}^m(1-q^{-2i})
\ge\frac34\left(1-\sum_{i=2}^\infty4^{-i}\right)
=\frac{11}{16}>\frac23.
\]
In orthogonal type, write $r=\lfloor n/2\rfloor$
and $d=\gcd(2,q-1)$.
The standard order formulas give
\[
\frac{zt}{q^D}=
\begin{cases}
\displaystyle
\frac1d\prod_{i=1}^r(1-q^{-2i}),
& n=2r+1,\\[6pt]
\displaystyle
\frac{1-\varepsilon q^{-r}}d
\prod_{i=1}^{r-1}(1-q^{-2i}),
& n=2r.
\end{cases}
\]
For every $k\ge1$,
\[
\prod_{i=1}^k(1-q^{-2i})
\ge1-\sum_{i=1}^kq^{-2i}
>1-\frac1{q^2-1}.
\]
For $n=2r$, we have $r\ge4$, so
$1-\varepsilon q^{-r}\ge1-q^{-4}$.
Since $d=2$ for odd $q$ and $d=1$
for even $q$, we have
\[
\frac{zt}{q^D}>
\begin{cases}
(1-3^{-4})(1-1/8)/2=35/81,
& q\text{ odd},\\
(1-2^{-4})(1-1/3)=5/8,
& q\text{ even}.
\end{cases}
\]
All these lower bounds exceed
$2/5>1/\delta$. Hence $\delta zt>q^D$.

We next bound $ot$.
By Table~\ref{tab:classical-parameters},
$o/z=2f$ in unitary type and $o/z=f$
in symplectic type.
In odd orthogonal type, $o=2f$ and $d=2$.
In even orthogonal type, $o/z=2fd$,
except for $O_8^+$, where $o/z=6fd$.
To bound $ot$, we use $ot=(o/z)zt$.
In unitary type, we have
$
(1-q^{-2i})(1+q^{-(2i+1)})<1,
$
and in orthogonal minus type,
$
(1-q^{-2})(1+q^{-r})<1.
$
Together with the order formulas above and
the values of $o/z$ in
Table~\ref{tab:classical-parameters},
these give
\[
ot<
\begin{cases}
fq^D,&\text{in odd orthogonal type},\\
6fq^D,&\text{for }O_8^+(q),\\
2fq^D,&\text{in all other cases}.
\end{cases}
\]
In particular, $ot<6fq^D$.
Finally, $z\le Q-1$ by
Table~\ref{tab:classical-parameters}.

\medskip
\noindent
\textit{Step 2. The higher-dimensional groups.}
By Lemma~\ref{lem:corefree}(a),
$
h\ge t\delta(Q-1)Q^{-\beta n}.
$
Combining this with $\delta zt>q^D$
from Step~1 gives
\[
\frac{h^{15/4}}{q^D}>
q^{(11D-15a\beta n)/4}
\left(\frac{Q-1}{z}\right)^{15/4}.
\]
Since $\beta<9/4$, we have
\[
\frac{h^{15/4}}{q^D}>
q^{A/16}\left(\frac{Q-1}{z}\right)^{15/4},
\qquad A=44D-135an.
\]

Suppose that $n\ge7$ in unitary type,
$n\ge6$ in symplectic type, or
$n\ge8$ in orthogonal type.
By Table~\ref{tab:classical-parameters},
\[
A=
\begin{cases}
44n^2-270n-44,&\text{unitary type},\\
22n^2-113n,&\text{symplectic type},\\
22n^2-157n,&\text{orthogonal type}.
\end{cases}
\]
These expressions are increasing in the
respective ranges, with initial values
$222$, $114$ and $152$.
In particular, $A\ge114$.
Since $z\le Q-1$, the bounds from Step~1
and $f\le q/2$ give
\[
\delta ot<18fq^D
\le9q^{D+1}
<q^{D+114/16}
\le q^{D+A/16}
<h^{15/4}.
\]
This proves the result in the stated ranges.

\medskip
\noindent
\textit{Step 3. The case $O_7(q)$.}
Here $q\ge3$ is odd, $Q=q$, $z=1$
and $A=-21$.
Since $q-1\ge2q/3$ and $2^8>3^5$,
we have $(q-1)^8>q^5$.
Thus the bound from Step~2 gives
\[
\frac{h^{15/4}}{q^D}
>q^{-21/16}(q-1)^{15/4}
>q^{33/32}
>q\ge3f
>\frac{\delta ot}{q^D},
\]
where the last inequality uses
$ot<fq^D$ and $\delta<3$.
Hence $h^{15/4}>\delta ot$.

\medskip
\noindent
\textit{Step 4. The cases $U_4(q)$ and $U_6(q)$.}
It remains to consider $U_4(q)$ with $q\ge3$
and $U_6(q)$ with $q\ge2$.
Write $q^\nu=|T|_p$.
The relevant parameters are as follows.
\begin{center}
\begin{minipage}{\linewidth}
\centering
\captionof{table}{Small unitary-group parameters}
\label{tab:unitary-small-parameters}
\[
\begin{array}{c|c|c|c|r|r|c}
T&\text{range}&\Phi_m(q)&j&a_0&\nu&c\\\hline
U_4(q)&q\ge3&q^2-q+1&3&16&6&1/4\\
U_6(q)&q\ge2&q^4-q^3+q^2-q+1&5&128&15&1/16
\end{array}
\]
\end{minipage}
\end{center}
Here $m=6,10$ in the respective rows,
$\Phi_m$ denotes the $m$-th cyclotomic
polynomial, and $j=m/2$.
The values of $\nu$ follow from
$|T|_p=q^{n(n-1)/2}$.
By Lemma~\ref{lem:corefree}(b), we have
$|U|_p\le p^{\lambda_p(n)}$ with $n=4,6$, respectively.
For $p=2$, Theorem~\ref{thm:dixonlin} gives
$\lambda_2(4)=4$ and $\lambda_2(6)=7$.
For odd $p$, the same formula gives
$p^{\lambda_p(4)}\le9$ and $p^{\lambda_p(6)}\le81$.
Hence $|U|_p\le a_0$ in both rows.
The constants $c$ are chosen for the
estimates below.
In these families $ot<2fq^D$, with
$(D,\nu)=(15,6),(35,15)$.
Since $\delta<3$ and $f\le q/2$, we have
$
\delta ot<3q^{D+1}.
$
Thus it suffices to prove
\begin{equation}\label{eq:cyclotomic-h}
h>cq^\nu
\end{equation}
and
$
c^{15}q^{15\nu-4D-4}>3^4.
$
Indeed, these inequalities imply
$h^{15/4}>3q^{D+1}>\delta ot$.

The exponents $15\nu-4D-4$ are
$26$ and $81$, respectively.
At the smallest allowed field sizes,
the second inequality becomes
$3^{22}>2^{30}$ and $2^{21}>3^4$,
respectively.
Since both exponents are positive,
it holds throughout the stated ranges.

To prove \eqref{eq:cyclotomic-h}, we choose
a prime $\ell\mid t$.
If $\ell\nmid|U|$, then $\ell\mid h$;
otherwise, we will use
Lemma~\ref{lem:corefree}(c) to bound $|U|$.
Write $Q=q^2$.
For $U_4(q)$, choose a prime $\ell\ne3$
dividing $\Phi_6(q)$.
For $U_6(q)$, choose a prime $\ell\ne5$
dividing $\Phi_{10}(q)$.
Such primes exist.
If $3\mid\Phi_6(q)$, then
$\Phi_6(q)\equiv3\pmod9$.
Similarly, if $5\mid\Phi_{10}(q)$,
then $\Phi_{10}(q)\equiv5\pmod{25}$.
Since $\Phi_6(q)>3$ for $q\ge3$
and $\Phi_{10}(q)>5$ for $q\ge2$,
the required primes exist.

Since $\Phi_m(q)$ is odd and
$\Phi_m(q)\equiv1\pmod p$, we have
$\ell\nmid2p$.
Our choice of $\ell$ also gives
$\ell\nmid m$.
Since $\ell\mid\Phi_m(q)$, the
cyclotomic-order criterion gives
\[
\operatorname{ord}_\ell(q)=m,
\qquad
\operatorname{ord}_\ell(Q)=j=m/2.
\]
In particular, $m\mid(\ell-1)$, so
$\ell\ge7$ for $U_4(q)$ and
$\ell\ge11$ for $U_6(q)$.
Thus $\ell>n+1$ and $n/2<j\le n$
in both cases.
Since $\Phi_6(q)\mid(q^3+1)$ and
$\Phi_{10}(q)\mid(q^5+1)$,
the order formula for $\SU_n(q)$ gives
$\ell\mid zt$.
As $\ell>n$ and $z=\gcd(n,q+1)$,
we have $\ell\mid t$.
We distinguish two cases.

\smallskip
\noindent
\textit{Case 1. $\ell\nmid|U|$.}
Then $\ell\mid h$, and
Lemma~\ref{lem:corefree}(b) gives
$h_p=q^\nu/|U|_p\ge q^\nu/a_0$.
Since $\ell\ne p$, we have
$h\ge\ell q^\nu/a_0$.
For $U_4(q)$ and $U_6(q)$, respectively,
$\ell/a_0\ge7/16>1/4=c$ and
$\ell/a_0\ge11/128>1/16=c$.
Hence \eqref{eq:cyclotomic-h} holds.

\smallskip
\noindent
\textit{Case 2. $\ell\mid|U|$.}
By Lemma~\ref{lem:corefree}(c),
\[
|U|\le
\frac{j(Q^j-1)|\GL_{n-j}(Q)|}{Q-1}
=
\begin{cases}
3(q^6-1),&T=U_4(q),\\
5(q^{10}-1),&T=U_6(q).
\end{cases}
\]

For $U_4(q)$, since $z\le q+1$ and
$(q^4-1)/(q^3-1)>q$, we have
\[
\frac{h}{q^6}
\ge\frac{(q^4-1)(q^2-1)}{3z(q^3-1)}
>\frac{q(q-1)}3>1.
\]

For $U_6(q)$, since $z\le q+1$ and
$(q^3+1)(q^2+1)>q^5-1$, we have
\[
\begin{aligned}
\frac{h}{q^{15}}
\ge
\frac{(q^2-1)(q^3+1)(q^4-1)(q^6-1)}
     {5z(q^5-1)}
>
\frac{(q-1)^2(q+1)(q^6-1)}5
\ge\frac{189}{5}>1.
\end{aligned}
\]
Thus $h>q^\nu>cq^\nu$ in both cases,
completing the proof.

\end{proof}

\begin{theorem}\label{thm:classical}
Let $T$ be one of the remaining classical
simple groups in Table~\ref{tab:classical-cases},
other than $\PSL_3(3)$.
Assume that \eqref{eq:inner-bound} holds
for every nonabelian simple group
of order less than $|T|$.
Then
\[
\delta|\Out(T)|\sigma(T)^{1+s}<|T|^s.
\]
\end{theorem}

\begin{proof}
Let $U\le T$ be soluble, and write $h=[T:U]>1$.
It suffices to prove
$
\delta|\Out(T)|\,|T|<h^{1+s}.
$
If $T=\PSL_n(q)$, this follows from
Proposition~\ref{prop:linear}.
If $T$ is one of $U_3(q)$, $U_5(q)$,
$S_4(q)$ or $U_4(2)$,
Lemma~\ref{lem:classindex} gives
$\delta|\Out(T)|\,|T|<h^{15/4}$.

For the remaining groups,
Lemma~\ref{lem:classPI} gives the required
inequality when $O_p(U)\ne1$.
When $O_p(U)=1$, Lemma~\ref{lem:classCF}
gives $\delta|\Out(T)|\,|T|<h^{15/4}$.
Since $15/4<1+s=53/14$, these bounds
establish the required inequality in
all cases.
Taking $|U|=\sigma(T)$ proves
\eqref{eq:outer-bound}.
\end{proof}

\section{Exceptional groups}\label{sec:exceptional}

We establish \eqref{eq:outer-bound} for the exceptional simple groups
of Lie type. For four families, the group with the smallest parameter
is not simple. Table~\ref{tab:exceptional-small} lists these cases and
gives the low-parameter identifications from \cite[Remark~24.18]{MT},
together with the earlier results that cover the resulting simple groups.

\begin{center}
\begin{minipage}{\linewidth}
\centering
\captionof{table}{Low-parameter cases of exceptional groups}
\label{tab:exceptional-small}
\[
\begin{array}{c|l|c}
\text{Low-parameter group}&\text{Relevant fact}&\text{Remaining range}\\\hline
G_2(2)
 &G_2(2)'\cong\PSU_3(3)
   \ \text{(Theorem~\ref{thm:classical})}
 &q\ge3\\
{}^2B_2(2)
 &\text{Soluble}
 &q\ge8\\
{}^2G_2(3)
 &{}^2G_2(3)'\cong\PSL_2(8)
   \ \text{(Proposition~\ref{prop:rankone})}
 &q\ge27\\
{}^2F_4(2)
 &{}^2F_4(2)'
   \ \text{(Proposition~\ref{prop:sporadic})}
 &q\ge8
\end{array}
\]
\end{minipage}
\end{center}

We first use minimum subgroup indices to handle seven families.
For the remaining families $E_6(q)$, $E_7(q)$ and $E_8(q)$,
we consider soluble subgroups $U\le T$ according to whether
$O_p(U)$ is trivial. The inductive hypothesis on
\eqref{eq:inner-bound} is needed only when $O_p(U)\ne1$.
Throughout this section, $q=p^f$ is the parameter in the group name.
For ${}^2B_2(q)$ and ${}^2F_4(q)$, we have $p=2$ and $f$ is odd.
For ${}^2G_2(q)$, we have $p=3$ and $f$ is odd.

\begin{lemma}\label{lem:exceptional-index}
Let $T$ be a simple group of one of the types
$G_2,F_4,{}^2E_6,{}^3D_4,{}^2F_4,{}^2B_2,{}^2G_2$,
with $q\ge3$ for $G_2(q)$,
$q\ge8$ for ${}^2B_2(q)$ and ${}^2F_4(q)$,
and $q\ge27$ for ${}^2G_2(q)$.
Then \eqref{eq:outer-bound} holds for $T$.
\end{lemma}

\begin{proof}
\medskip
\noindent
\textit{Step 1. Reduction to minimum indices.}
Write $o=|\Out(T)|$ and
$
\mu(T)=\min\{[T:H]:H<T\}.
$
Since every soluble subgroup of $T$ is proper,
$\sigma(T)\le |T|/\mu(T)$.
Thus, to prove \eqref{eq:outer-bound}, it suffices to show
$
\delta o|T|<\mu(T)^{53/14}.
$
We first consider all cases except $G_2(4)$
and ${}^2E_6(q)$, which will be handled
separately in Step~4.

\medskip
\noindent
\textit{Step 2. Bounds on $\mu(T)$ and $o|T|$.}
For the cases treated in Step~3, we establish
bounds of the form
\[
\mu(T)>q^m,\qquad o|T|<cfq^D,
\]
where $c,D,m$ depend on the family.

For simple $T$, the minimum faithful permutation degree
equals $\mu(T)$, as noted in the proof of
Lemma~\ref{lem:classindex}.
The required minimum indices are given by
Vasil'ev~\cite{VasGF}
for $G_2$ and $F_4$, and by
Vasil'ev~\cite{VasTw}
for the five twisted families; they are also listed in
\cite[Table~4]{GMPS}:
\begin{center}
\begin{minipage}{\linewidth}
\centering
\captionof{table}{Minimum indices of exceptional groups}
\label{tab:exceptional-min-indices}
\[
\begin{array}{c|l}
 T&\mu(T)\\\hline
 G_2(3),\ G_2(4)&351,\ 416\\
 G_2(q),\ q>4&(q^6-1)/(q-1)\\
 F_4(q)&(q^{12}-1)(q^4+1)/(q-1)\\
 {}^2E_6(q)&(q^{12}-1)(q^6-q^3+1)(q^4+1)/(q-1)\\
 {}^3D_4(q)&(q^8+q^4+1)(q+1)\\
 {}^2F_4(q)&(q^6+1)(q^3+1)(q+1)\\
 {}^2B_2(q)&q^2+1\\
 {}^2G_2(q)&q^3+1
\end{array}
\]
\end{minipage}
\end{center}
These formulas give
\[
\begin{aligned}
 \mu(G_2(q))&>q^5 &&(q\ne4),\\
 \mu(F_4(q))&>q^{15},\\
 \mu({}^3D_4(q))&>q^9,\\
 \mu({}^2F_4(q))&>q^{10},\\
 \mu({}^2B_2(q))&>q^2,\\
 \mu({}^2G_2(q))&>q^3.
\end{aligned}
\]
For $G_2(3)$, the bound follows from
$351>3^5$. It fails for $G_2(4)$, since
$\mu(G_2(4))=416<4^5$.
This case will be treated separately in Step~4.

We next establish the upper bounds on $o|T|$,
using the standard group-order and
outer-automorphism formulas~\cite{MT,Wilson}.
For $G_2(q)$ and $F_4(q)$, respectively,
\[
\begin{aligned}
 |T|&=q^6(q^6-1)(q^2-1),\\
 |T|&=q^{24}(q^{12}-1)(q^8-1)(q^6-1)(q^2-1).
\end{aligned}
\]
Since $o\le2f$ in both cases, we obtain
\[
o|G_2(q)|<2fq^{14},\qquad
o|F_4(q)|<2fq^{52}.
\]
For the four twisted families under consideration,
the orders are
\begin{center}
\begin{minipage}{\linewidth}
\centering
\captionof{table}{Orders of twisted exceptional groups}
\label{tab:twisted-exceptional-orders}
\[
\begin{array}{c|l|c}
 T&|T|&o\\\hline
 {}^3D_4(q)&q^{12}(q^8+q^4+1)(q^6-1)(q^2-1)&3f\\
 {}^2F_4(q)&q^{12}(q^6+1)(q^4-1)(q^3+1)(q-1)&f\\
 {}^2B_2(q)&q^2(q^2+1)(q-1)&f\\
 {}^2G_2(q)&q^3(q^3+1)(q-1)&f
\end{array}
\]
\end{minipage}
\end{center}
For ${}^3D_4(q)$, we have
\[
(1-q^{-2})(1+q^{-4}+q^{-8})<1.
\]
For ${}^2F_4(q)$, we have
\[
(1-q^{-4})(1+q^{-6})<1,\qquad
(1-q^{-1})(1+q^{-3})<1.
\]
For ${}^2B_2(q)$ and ${}^2G_2(q)$,
the corresponding inequalities are
\[
(1-q^{-1})(1+q^{-2})<1,\qquad
(1-q^{-1})(1+q^{-3})<1,
\]
respectively.
The preceding formulas and inequalities give
the following bounds for the six families
considered in Step~3:
\begin{center}
\begin{minipage}{\linewidth}
\centering
\captionof{table}{Bounds for exceptional groups}
\label{tab:exceptional-index-bounds}
\[
\begin{array}{c|c|c}
 T&\mu(T)&o|T|\\\hline
 G_2(q),\ q\ne4&>q^5&<2fq^{14}\\
 F_4(q)&>q^{15}&<2fq^{52}\\
 {}^3D_4(q)&>q^9&<3fq^{28}\\
 {}^2F_4(q)&>q^{10}&<fq^{26}\\
 {}^2B_2(q)&>q^2&<fq^5\\
 {}^2G_2(q)&>q^3&<fq^7
\end{array}
\]
\end{minipage}
\end{center}

\medskip
\noindent
\textit{Step 3. Verification of the required inequality.}
The bounds established in Step~2 give the following
values of $c,D,m$ and $M=53m-14D$:
\begin{center}
\begin{minipage}{\linewidth}
\centering
\captionof{table}{Parameters for the exponent comparisons}
\label{tab:exceptional-index-parameters}
\[
\begin{array}{c|r|r|r|r}
 T&c&D&m&M\\\hline
 G_2(q),\ q\ne4&2&14&5&69\\
 F_4(q)&2&52&15&67\\
 {}^3D_4(q)&3&28&9&85\\
 {}^2F_4(q)&1&26&10&166\\
 {}^2B_2(q)&1&5&2&36\\
 {}^2G_2(q)&1&7&3&61
\end{array}
\]
\end{minipage}
\end{center}
Since $q\ge2$, the values in
Table~\ref{tab:exceptional-index-parameters} give
$
(3c/2)^{14}<2^{M-14}\le q^{M-14}.
$
Together with $f\le q/2$ and $\delta<3$, this yields
\[
\begin{aligned}
(\delta o|T|)^{14}
<(3cfq^D)^{14}
\le(3c/2)^{14}q^{14D+14}
<q^{53m}
<\mu(T)^{53}.
\end{aligned}
\]
Thus $\delta o|T|<\mu(T)^{53/14}$
for all groups in the table.

\medskip
\noindent
\textit{Step 4. The remaining cases.}
It remains to consider $G_2(4)$ and ${}^2E_6(q)$.

\smallskip
\noindent
\textit{The case $G_2(4)$.}
The values
\[
|T|=251596800,\qquad o=2,\qquad \mu(T)=416
\]
give
$
(6\cdot251596800)^{14}<416^{53}.
$
Since $\delta<3$, this proves
$\delta o|T|<\mu(T)^{53/14}$.

\smallskip
\noindent
\textit{The family ${}^2E_6(q)$.}
Write $d=\gcd(3,q+1)$.
The group order and outer-automorphism order satisfy
\[
\begin{aligned}
 d|T|=q^{36}(q^2-1)(q^5+1)(q^6-1)
                  (q^8-1)(q^9+1)(q^{12}-1),\qquad
 o=2fd.
\end{aligned}
\]
Since
\[
(1-q^{-2})(1+q^{-5})<1,\qquad
(1-q^{-8})(1+q^{-9})<1,
\]
we have $d|T|<q^{78}$. Hence
$
o|T|=2f(d|T|)<2fq^{78}\le q^{79}.
$
The minimum-index formula gives
\[
\begin{aligned}
 \mu(T)>(q^{11}+q^{10}+q^9)(q^6-q^3)q^4
 =q^{20}(q+1)+q^{16}(q^3-q^2-q-1)
 >q^{20}(q+1),
\end{aligned}
\]
since
$
q^3-q^2-q-1
=(q-2)(q^2+q+1)+1>0
$
for $q\ge2$.
As $53/14>3$, we have
\[
\begin{aligned}
 \mu(T)^{53/14}
 >q^{79}\sqrt q(1+1/q)^{53/14}
 >q^{79}\left(\sqrt q+\frac3{\sqrt q}\right)
 \ge2\sqrt3\,q^{79}
 >\delta o|T|.
\end{aligned}
\]
Thus $\delta o|T|<\mu(T)^{53/14}$
also in these two cases, and
\eqref{eq:outer-bound} follows.
\end{proof}

\begin{lemma}\label{lem:exceptional-split}
Let $T=E_6(q)$, $E_7(q)$ or $E_8(q)$.
Suppose that \eqref{eq:inner-bound} holds for every
nonabelian simple group of order less than $|T|$.
Then \eqref{eq:outer-bound} holds for $T$.
\end{lemma}

\begin{proof}
Write $o=|\Out(T)|$.
Since $\sigma(T)$ is the maximum order of a soluble
subgroup of $T$, to prove \eqref{eq:outer-bound}
it suffices to show that
$
\delta o|U|^{1+s}<|T|^s
$
for every soluble subgroup $U\le T$.
Let $U\le T$ be an arbitrary soluble subgroup.
As in Section~\ref{sec:lie}, let $X=\mathbf G^\sigma$
be the simply connected cover of $T$, and write
$z=|Z(X)|$.

\medskip
\noindent
\textit{Case 1. $O_p(U)\ne1$.}

\smallskip
\noindent
\textit{Step 1. Reduction to a parabolic estimate.}
By Theorem~\ref{thm:bt}, the full preimage $\widehat U$
lies in a maximal parabolic subgroup $\widehat P$ of $X$.
Choose a Levi decomposition $\widehat P=QL$, where
$Q$ is the unipotent radical and $L=\mathbf L^\sigma$ is a
Levi complement, with $\mathbf L$ a $\sigma$-stable Levi
subgroup of $\mathbf G$.
Write $|Q|=q^u$.
By Proposition~\ref{prop:parabolic} and \eqref{eq:PIprime},
it is enough to establish
\[
Q=O_p(\widehat P)
\qquad\text{and}\qquad
3(o/z)c(L)\le q^{25u/14}.
\]

\smallskip
\noindent
\textit{Step 2. Levi structure and unipotent exponents.}
We use the standard structure of maximal parabolic
subgroups of split groups of Lie type
~\cite[Chapter~12]{MT}.
By~\cite[Proposition~12.14]{MT},
the derived subgroup $[\mathbf L,\mathbf L]$
is simply connected.
Hence $K=[\mathbf L,\mathbf L]^\sigma$
is a direct product of simply connected groups.
Their types correspond to the components of
the Dynkin diagram obtained by deleting one node.

The quotient $\mathbf L/[\mathbf L,\mathbf L]$
is a split torus of rank one.
Since $[\mathbf L,\mathbf L]$ is connected,
the Lang--Steinberg theorem
~\cite[Theorem~21.7]{MT}
gives a surjection
$
L\to(\mathbf L/[\mathbf L,\mathbf L])^\sigma
$
with kernel $K$.
Thus
\[
|L/K|=q-1,\qquad c(L)=(q-1)c(K).
\]
Every direct factor of $K$ has trivial $p$-core.
This includes $\SL_2(2)$ and $\SL_2(3)$,
which may occur for small $q$.
Hence $O_p(K)=1$.
Since $p\nmid(q-1)$, we have $O_p(L)=1$.
Consequently, $Q=O_p(\widehat P)$.

The parabolic subgroup $\widehat P$ contains
a Sylow $p$-subgroup of $X$.
Since $|L/K|=q-1$, we have
\[
q^u=\frac{|X|_p}{|K|_p}.
\]
For a split simply connected group,
the $p$-part of its order is $q^N$,
where $N$ is the number of positive roots.
Thus $u$ is the difference between the numbers
of positive roots of $\mathbf G$ and
$[\mathbf L,\mathbf L]$.

In Bourbaki numbering, the nodes
$1,3,4,5,\ldots,\ell$ form a chain,
and node $2$ is attached to node $4$.
The types of $K$ and the corresponding values of $u$
are given in Table~\ref{tab:exceptional-levi}
~\cite{Carter,MT,GLS3,Wilson}.
\begin{center}
\begin{minipage}{\linewidth}
\centering
\captionof{table}{Levi types and values of $u$ for split exceptional groups}
\label{tab:exceptional-levi}
\[
\begin{array}{c|lr|lr|lr}
 &\multicolumn{2}{c|}{E_6(q)}
 &\multicolumn{2}{c|}{E_7(q)}
 &\multicolumn{2}{c}{E_8(q)}\\
 \text{Deleted node}
 &\text{Type of }K&u
 &\text{Type of }K&u
 &\text{Type of }K&u\\\hline
 1&D_5&16&D_6&33&D_7&78\\
 2&A_5&21&A_6&42&A_7&92\\
 3&A_1A_4&25&A_1A_5&47&A_1A_6&98\\
 4&A_2A_1A_2&29&A_2A_1A_3&53&A_2A_1A_4&106\\
 5&A_4A_1&25&A_4A_2&50&A_4A_3&104\\
 6&D_5&16&D_5A_1&42&D_5A_2&97\\
 7&-&-&E_6&27&E_6A_1&83\\
 8&-&-&-&-&E_7&57
\end{array}
\]
\end{minipage}
\end{center}
In particular, $u\ge16$ in every case.

\smallskip
\noindent
\textit{Step 3. Verification of the parabolic estimate.}
Let $J$ be a direct factor of $K$ in the above decomposition.
By \cite[Theorem~24.17]{MT}, the only factors in the table that are
not quasisimple are $\SL_2(2)$ and $\SL_2(3)$. Both are soluble and
have order at most $24$. For every other factor $J$, we have
$c(J)=|Z(J)|\le24$ by \cite[Table~24.2]{MT}.
Since $K$ has at most three such factors,
$c(K)\le24^3$. Hence
\[
c(L)=(q-1)c(K)\le(q-1)24^3.
\]
The outer-automorphism orders~\cite{MT,Wilson}
give $o/z\le2f$.
Since $f\le q/2$, $u\ge16$ and
$3\cdot24^3<2^{16}\le q^{16}$, we have
\[
\begin{aligned}
3(o/z)c(L)
\le6f(q-1)24^3
<3\cdot24^3q^2
<q^{18}
<q^{25u/14}.
\end{aligned}
\]
Thus \eqref{eq:PIprime} holds.
Applying Proposition~\ref{prop:parabolic}
with the assumed inner inequalities gives
$\delta o|U|^{1+s}<|T|^s$.

\medskip
\noindent
\textit{Case 2. $O_p(U)=1$.}
By Lemma~\ref{lem:Ecorefree}, we have $|U|\le q^a$,
where $a=52,91,104$ for $E_6(q),E_7(q),E_8(q)$,
respectively.
We establish a lower bound on $|T|$ and an upper bound
on $\delta o$, and then prove the required inequality.

\smallskip
\noindent
\textit{Step 1. Bounds on $|T|$ and $\delta o$.}
The standard order formulas~\cite{MT,Wilson} give
$
z|T|=q^D\prod_{i\in I}(1-q^{-i}),
$
with the following data:
\begin{center}
\begin{minipage}{\linewidth}
\centering
\captionof{table}{Order parameters for split exceptional groups}
\label{tab:split-exceptional-orders}
\[
\begin{array}{c|r|c|l|c}
 T&D&z&I&o/z\\\hline
 E_6(q)&78&\gcd(3,q-1)&2,5,6,8,9,12&2f\\
 E_7(q)&133&\gcd(2,q-1)&2,6,8,10,12,14,18&f\\
 E_8(q)&248&1&2,8,12,14,18,20,24,30&f
\end{array}
\]
\end{minipage}
\end{center}
Since $I$ consists of distinct integers at least two,
\[
\prod_{i\in I}(1-q^{-i})
>1-\sum_{i\in I}q^{-i}
>1-\frac1{q(q-1)}
\ge\frac12.
\]
Since $z\le3$ and $q\ge2$, it follows that
$
|T|>q^D/(2z)\ge q^D/6>q^{D-3}.
$
The outer-automorphism orders in the table give
$o/z\le2f$. Since $z\le3$, $\delta<3$ and
$f\le q/2$, we have
$
\delta o<18f\le9q<q^5.
$

\smallskip
\noindent
\textit{Step 2. Verification of the required inequality.}
The values of $a$ and $D$ give the following
exponent comparisons:
\begin{center}
\begin{minipage}{\linewidth}
\centering
\captionof{table}{Exponent comparisons for split exceptional groups}
\label{tab:split-exceptional-comparisons}
\[
\begin{array}{c|r|r|rcl}
 T&a&D&53a+70&&39(D-3)\\\hline
 E_6(q)&52&78&2826&<&2925\\
 E_7(q)&91&133&4893&<&5070\\
 E_8(q)&104&248&5582&<&9555
\end{array}
\]
\end{minipage}
\end{center}
Combining these inequalities with the bounds on
$|U|$, $\delta o$ and $|T|$, we have
\[
(\delta o)^{14}|U|^{53}
<q^{70+53a}
<q^{39(D-3)}
<|T|^{39}.
\]
Thus $\delta o|U|^{1+s}<|T|^s$ also in this case.
\end{proof}

\section{Completion of the proof}\label{sec:proofs}

\begin{proof}[Proof of Theorem~\ref{thm:uniform}]
We prove \eqref{eq:inner-bound} by strong induction on $|T|$
and determine at each order whether \eqref{eq:outer-bound} holds.
Assume that \eqref{eq:inner-bound} holds for every
nonabelian simple group of order less than $|T|$.
By the argument in Proposition~\ref{prop:parabolic},
the nonabelian composition factors of the Levi quotient
have order less than $|T|$.

By the induction hypothesis, we may apply
Theorem~\ref{thm:classical} and
Lemma~\ref{lem:exceptional-split}.
Together with Propositions~\ref{prop:rankone},
\ref{prop:alt} and~\ref{prop:sporadic}
and Lemma~\ref{lem:exceptional-index},
these results cover all nonabelian finite simple groups,
by the classification of finite simple groups and the
reductions in Sections~\ref{sec:classical}
and~\ref{sec:exceptional}.
They establish \eqref{eq:outer-bound} unless
$T\cong\PSL_2(8)$, $\PSL_2(16)$ or $\PSL_3(3)$.
Since $\delta|\Out(T)|>1$,
\eqref{eq:outer-bound} implies \eqref{eq:inner-bound}.

It remains to consider the three exceptions.
The values of $\sigma(T)$ for $\PSL_2(8)$ and
$\PSL_2(16)$ were obtained in the proof of
Proposition~\ref{prop:rankone},
and that for $\PSL_3(3)$ follows from
\eqref{eq:mann-sigma}.
The relevant values are as follows:
\begin{center}
\begin{minipage}{\linewidth}
\centering
\captionof{table}{The three exceptions to \eqref{eq:outer-bound}}
\label{tab:outer-bound-exceptions}
\[
\begin{array}{c|rrr}
 T&\sigma(T)&|T|/\sigma(T)&|\Out(T)|\\\hline
 \PSL_2(8)&56&9&3\\
 \PSL_2(16)&240&17&4\\
 \PSL_3(3)&432&13&2
\end{array}
\]
\end{minipage}
\end{center}
For each row of the table, a direct calculation gives
\[
\sigma(T)<(|T|/\sigma(T))^s
<\delta|\Out(T)|\sigma(T).
\]
The first inequality is equivalent to \eqref{eq:inner-bound},
and the second shows that \eqref{eq:outer-bound} fails.
This completes the induction and the proof of both assertions
of the theorem.
\end{proof}

\begin{proof}[Proof of Theorem~A]
Let $T$ be a nonabelian finite simple group and let
$B\le\Aut(T)$ be soluble. Write $h=[T:B\cap T]$.
Since $B$ is soluble and $T$ is nonsoluble, we have $h>1$.

If $T$ is not one of the three exceptions in
Theorem~\ref{thm:uniform}, then
\eqref{eq:outer-bound} and Lemma~\ref{lem:crude}
give $\delta|B|<h^s<h^{r_0}$.
For the three exceptions,
Remarks~\ref{rem:rankone-exceptions}
and~\ref{rem:L33-exception} give
$\delta|B|\le h^{r_0}$.
Thus
\[
\delta|B|\le[T:B\cap T]^{r_0}
\]
for every nonabelian finite simple group $T$
and every soluble $B\le\Aut(T)$.
Since $b_0=1+r_0$, Theorem~B with $r=r_0$
gives both inequalities in Theorem~A.

Now let $b$ be a real number such that
$[X:\Rad(X)]\le[X:H]^b$ for every finite group $X$
and every soluble subgroup $H\le X$.
Taking $X$ to be a nonabelian finite simple group
and $H=1$ gives $b\ge1$.
By Theorem~B with $r=b-1$, we have
\[
\delta|B|\le[T:B\cap T]^{b-1}
\]
for every nonabelian finite simple group $T$
and every soluble $B\le\Aut(T)$.
Taking $T$ and $B$ as in Proposition~\ref{prop:lower},
we have
$
168\delta\le9^{b-1}.
$
Hence $b\ge1+\log_9(168\delta)=b_0$.

Finally, the equality attained in
Proposition~\ref{prop:lower} shows that,
for the fixed exponent $b_0$, the coefficient
$\delta^{-1}$ in the second inequality of
Theorem~A cannot be decreased.
\end{proof}

\appendix
\section{Finite arithmetic verification}\label{app:verification}
We give the remaining arithmetic checks and justify the list of
classical parabolic cases.

\subsection{The groups \texorpdfstring{$\PSL_2(q)$}{PSL(2,q)}}\label{app:rankone}

In Proposition~\ref{prop:rankone}, it remains to verify
$\delta fq(q^2-1)<k^{1+s}$ for the prime powers
$q=p^f$ with $4\le q<64$ and $q\notin\{5,8,16\}$.
The corresponding values of $k$ and $f$ and the
required inequalities are given below.
\begin{center}
\begin{minipage}{\linewidth}
\centering
\captionof{table}{Arithmetic checks for $\PSL_2(q)$}
\label{tab:rankone-arithmetic}
\[
\begin{array}{c|c|c|l}
 q&k&f&\delta fq(q^2-1)<k^{1+s}\\\hline
 4&5&2&120\delta<5^{1+s}\\
 7&7&1&336\delta<7^{1+s}\\
 11\le q<64,\ q\text{ prime}&q+1&1&
 \delta q(q^2-1)<(q+1)^{1+s}\\
 9,25,49&q+1&2&
 2\delta q(q^2-1)<(q+1)^{1+s}\\
 27&q+1&3&
 3\delta q(q^2-1)<(q+1)^{1+s}\\
 32&q+1&5&
 5\delta q(q^2-1)<(q+1)^{1+s}
\end{array}
\]
\end{minipage}
\end{center}

The inequalities in the first two rows hold
by direct calculation.
In the remaining four rows, since $\delta<3$
and $q(q-1)<(q+1)^2$, we have
$\delta fq(q^2-1)<3f(q+1)^3$.
Since $s-2=11/14$, the required inequality
follows from $(3f)^{14}<(q+1)^{11}$.
A direct calculation shows that this inequality holds
at the smallest value of $q$ in each row,
and hence for all the listed values of $q$.

\subsection{Remaining classical parabolic cases}
\label{app:classical-parabolic}
In Lemma~\ref{lem:classPI}, the remaining cases satisfy
$u\le8$ when $q=2$ and $u\le7$ when $q=3$.
Since $u\ge n-2$, we obtain the corresponding bounds
on $n$. For $U_n(q)$ with $n\ge6$, we have $u\ge9$.
Together with the restrictions in Lemma~\ref{lem:classPI},
these bounds give the following possible values of $n$:
\begin{center}
\begin{minipage}{\linewidth}
\centering
\captionof{table}{Possible dimensions for the remaining cases}
\label{tab:classical-parabolic-dimensions}
\[
\begin{array}{c|cc|cccc}
 &\multicolumn{2}{c|}{\text{Bounds}}
 &\multicolumn{4}{c}{\text{Possible values of }n}\\
 q&u&n&U_n&S_n&O_n\ (n\text{ odd})&O_n^\pm\\\hline
 2&\le8&\le10&-&6,8,10&-&8,10\\
 3&\le7&\le9&4&6,8&7,9&8
\end{array}
\]
\end{minipage}
\end{center}

For each group in Table~\ref{tab:classical-parabolic-dimensions},
we use the ranges of $k$ and the formulas for $u$ in
Table~\ref{tab:classical-levi}.
The resulting values of $u$ are given below.
The last two columns list the values of $k$ that
remain to be considered for $q=2$ and $q=3$, respectively.

\begin{center}
\begin{minipage}{\linewidth}
\centering
\captionof{table}{Enumeration of the remaining classical parabolic cases}
\label{tab:classical-parabolic-enumeration}
\[
\begin{array}{c|l|l|cc}
 &&&\multicolumn{2}{c}{\text{Remaining values of }k}\\
 \text{Type}&k&u&q=2&q=3\\\hline
 U_4&1,2&5,4&-&1,2\\
 S_6&1,2,3&5,7,6&1,2,3&1,2,3\\
 S_8&1,2,3,4&7,11,12,10&1&1\\
 S_{10}&1,2,3,4,5&9,15,18,18,15&-&-\\
 O_7&1,2,3&5,7,6&-&1,2,3\\
 O_9&1,2,3,4&7,11,12,10&-&1\\
 O_8^+&1,2,4&6,9,6&1,4&1,4\\
 O_8^-&1,2,3&6,9,9&1&1\\
 O_{10}^+&1,2,3,5&8,13,15,10&1&-\\
 O_{10}^-&1,2,3,4&8,13,15,14&1&-
\end{array}
\]
\end{minipage}
\end{center}

For type $O_n^+$, the two classes of maximal totally
singular subspaces are counted once because they give
the same values of $u$ and $c(L)$.
This gives the nine cases for $q=2$ and the
thirteen cases for $q=3$ listed in
Table~\ref{tab:classical-parabolic-cases}.

\section*{Declaration of AI use}

The author used ChatGPT (GPT-6, OpenAI) and Claude
(Opus 5.5, Anthropic) to assist with language editing,
discussions of mathematical arguments, and some of
the more involved calculations and estimates for
finite groups of Lie type.
All mathematical results and proofs in the final
manuscript were verified by the author, who takes
full responsibility for its content.

\bibliographystyle{amsplain}
\bibliography{BGP_references}

@article{BGP,
  author  = {Babai, L. and Goodman, A. J. and Pyber, L.},
  title   = {Groups without faithful transitive permutation representations of small degree},
  journal = {J. Algebra},
  volume  = {195},
  year    = {1997},
  pages   = {1--29},
}

@article{BW,
  author  = {Burgoyne, N. and Williamson, C.},
  title   = {On a theorem of {Borel} and {Tits} for finite {Chevalley} groups},
  journal = {Arch. Math. (Basel)},
  volume  = {27},
  year    = {1976},
  pages   = {489--491},
}

@article{Burness,
  author  = {Burness, T. C.},
  title   = {On soluble subgroups of sporadic groups},
  journal = {Israel J. Math.},
  volume  = {254},
  year    = {2023},
  pages   = {313--340},
}

@article{BurnessPrimitive,
  author  = {Burness, T. C.},
  title   = {Base sizes for primitive groups with soluble stabilisers},
  journal = {Algebra Number Theory},
  volume  = {15},
  year    = {2021},
  pages   = {1755--1807},
}

@article{BastosSchneider,
  author  = {Bastos, R. and Schneider, C.},
  title   = {Groups with a solvable subgroup of prime-power index},
  journal = {Monatsh. Math.},
  volume  = {194},
  year    = {2021},
  pages   = {471--480},
}

@book{Carter,
  author    = {Carter, R. W.},
  title     = {Finite groups of {Lie} type: conjugacy classes and complex characters},
  publisher = {Wiley},
  year      = {1985},
}

@article{Dixon,
  author  = {Dixon, J. D.},
  title   = {The {Fitting} subgroup of a linear solvable group},
  journal = {J. Austral. Math. Soc.},
  volume  = {7},
  year    = {1967},
  pages   = {417--424},
}

@article{DixonNormal,
  author  = {Dixon, J. D.},
  title   = {Normal {$p$}-subgroups of solvable linear groups},
  journal = {J. Austral. Math. Soc.},
  volume  = {7},
  year    = {1967},
  pages   = {545--551},
}

@article{GMPS,
  author  = {Guest, S. and Morris, J. and Praeger, C. E. and Spiga, P.},
  title   = {On the maximum orders of elements of finite almost simple groups and primitive permutation groups},
  journal = {Trans. Amer. Math. Soc.},
  volume  = {367},
  year    = {2015},
  pages   = {7665--7694},
}

@book{GLS3,
  author    = {Gorenstein, D. and Lyons, R. and Solomon, R.},
  title     = {The classification of the finite simple groups, number~3},
  series    = {Math. Surveys Monogr.},
  volume    = {40.3},
  publisher = {Amer. Math. Soc.},
  year      = {1998},
}

@book{Huppert,
  author    = {Huppert, B.},
  title     = {Endliche {Gruppen} {I}},
  publisher = {Springer},
  year      = {1967},
}

@book{Kourovka,
  editor    = {Khukhro, E. I. and Mazurov, V. D.},
  title     = {Unsolved problems in group theory: The {Kourovka Notebook}, no.~21},
  publisher = {Sobolev Institute of Mathematics},
  address   = {Novosibirsk},
  year      = {2026},
  note      = {arXiv:1401.0300},
}

@book{KL,
  author    = {Kleidman, P. B. and Liebeck, M. W.},
  title     = {The subgroup structure of the finite classical groups},
  series    = {London Math. Soc. Lecture Note Ser.},
  volume    = {129},
  publisher = {Cambridge Univ. Press},
  year      = {1990},
}

@article{Mann,
  author  = {Mann, A.},
  title   = {Soluble subgroups of symmetric and linear groups},
  journal = {Israel J. Math.},
  volume  = {55},
  year    = {1986},
  pages   = {162--172},
}

@book{MT,
  author    = {Malle, G. and Testerman, D.},
  title     = {Linear algebraic groups and finite groups of {Lie} type},
  series    = {Cambridge Stud. Adv. Math.},
  volume    = {133},
  publisher = {Cambridge Univ. Press},
  year      = {2011},
}

@article{Palfy,
  author  = {P{\'a}lfy, P. P.},
  title   = {A polynomial bound for the orders of primitive solvable groups},
  journal = {J. Algebra},
  volume  = {77},
  year    = {1982},
  pages   = {127--137},
}

@article{VasGF,
  author  = {Vasil'ev, A. V.},
  title   = {Minimal permutation representations of finite simple exceptional groups of types {$G_2$} and {$F_4$}},
  journal = {Algebra Logic},
  volume  = {35},
  year    = {1996},
  pages   = {371--383},
}

@article{VasTw,
  author  = {Vasil'ev, A. V.},
  title   = {Minimal permutation representations of finite simple exceptional twisted groups},
  journal = {Algebra Logic},
  volume  = {37},
  year    = {1998},
  pages   = {9--20},
}

@article{Vdovin,
  author  = {Vdovin, E. P.},
  title   = {The number of subgroups with trivial unipotent radicals in finite groups of {Lie} type},
  journal = {J. Group Theory},
  volume  = {7},
  year    = {2004},
  pages   = {99--112},
}

@book{Wilson,
  author    = {Wilson, R. A.},
  title     = {The finite simple groups},
  series    = {Grad. Texts in Math.},
  volume    = {251},
  publisher = {Springer},
  year      = {2009},
}

@article{Wolf,
  author  = {Wolf, T. R.},
  title   = {Solvable and nilpotent subgroups of {$\mathrm{GL}(n,q^m)$}},
  journal = {Canad. J. Math.},
  volume  = {34},
  year    = {1982},
  pages   = {1097--1111},
}

@book{ATLAS,
  author    = {Conway, J. H. and Curtis, R. T. and Norton, S. P. and Parker, R. A. and Wilson, R. A.},
  title     = {Atlas of finite groups: maximal subgroups and ordinary characters for simple groups},
  publisher = {Clarendon Press},
  address   = {Oxford},
  year      = {1985},
  isbn      = {0-19-853199-0},
}
\end{document}